\documentclass[10pt]{amsart}
\usepackage{rotating}
\usepackage{mathcomp,amscd}
\usepackage{amssymb}
\usepackage{mathtools}
\usepackage[all]{xy}
\usepackage{subfigure}
\usepackage[margin=1cm]{caption}
\usepackage{wasysym}
\usepackage{color}
\definecolor{dg}{rgb}{0,0.67,0}
\definecolor{gr}{rgb}{.67,0.67,0.67}
\usepackage{amsfonts}
\usepackage{url}
\usepackage{graphicx}
\usepackage{amssymb}
\usepackage{epstopdf}
\usepackage{array}

  \makeindex

\newcommand{\Out}{\mathrm{Out}}
\newcommand{\Spec}{\mathrm{Spec}}

\newcommand{\Hur}{{\mbox{\sc Hur}}}

\newcommand{\AHur}{{\sf Hur}}

\renewcommand{\AC}{{\sf C}}
\newcommand{\Conf}{{\mbox{\sc Conf}}}
\newcommand{\scX}{{\mbox{\sc X}}}
\newcommand{\scU}{{\mbox{\sc U}}}
\newcommand{\AConf}{{\sf Conf}}
\newcommand{\AQuart}{{\sf Quart}}

\newcommand{\TopcircleA}{\mbox{\begin{sideways} {\!\Rightcircle \,} \end{sideways}}}

\usepackage[OT2,T1]{fontenc}
\DeclareSymbolFont{cyrletters}{OT2}{wncyr}{m}{n}
\DeclareFontFamily{OT1}{rsfs}{}
\DeclareFontEncoding{OT2}{}{} 
  
     \DeclareFontShape{OT1}{rsfs}{n}{it}{<-> rsfs10}{}
\DeclareMathAlphabet{\mathscr}{OT1}{rsfs}{n}{it}

\newcommand{\C}{{\Bbb C}}
\newcommand{\Z}{{\Bbb Z}}
\newcommand{\R}{{\Bbb R}}

\newcommand{\cG}{\mathcal{G}}

\newcommand{\F}{{\Bbb F}}

\newcommand{\bbP}{{\Bbb P}}
\newcommand{\AP}{{\sf P}}  

\newcommand{\AX}{{\sf X}}

\newcommand{\SB}{{\mbox{\sc B}}}
\newcommand{\SC}{{\mbox{\sc C}}}

\newcommand{\AAnew}{{\mbox{\sf A}}}
\newcommand{\AB}{{\mbox{\sf B}}}
\newcommand{\AD}{{\mbox{\sf D}}}

\newcommand{\AS}{{\sf S}}

\newcommand{\AU}{{\sf U}}
\newcommand{\SU}{{\mbox{\sc U}}}

\newcommand{\Q}{{\Bbb Q}}
 
\newcommand{\Gal}{\mbox{\rm Gal}}
\newcommand{\Aut}{\mbox{Aut}}

\newcommand{\cP}{\mathcal{P}}
\newcommand{\APGL}{{\sf PGL}}

\newcommand{\disc}{\mbox{disc}}

\newcommand{\ord}{{\rm ord}}

\numberwithin{equation}{section}
\numberwithin{table}{section}
\numberwithin{figure}{section}
\newtheorem{Theorem}{Theorem}[section]
\newtheorem*{Theorem*}{Theorem}
\newtheorem*{corollary*}{Corollary}

\newtheorem*{Proposition*}{Proposition}

\newtheorem{Conjecture}[Theorem]{Conjecture}

\newcommand{\cmmt}[1]{}
 \allowdisplaybreaks
\title{Hurwitz Number Fields}
\author{David P. Roberts} 
\address{Division of Science and Mathematics, University of Minnesota Morris; 
Morris, Minnesota, 56267, USA}
\email{roberts@morris.umn.edu}
\urladdr{http://cda.morris.umn.edu/~roberts}
\begin{document}

\begin{abstract}  
The canonical covering maps from Hurwitz varieties to configuration varieties
are important in algebraic geometry.
The scheme-theoretic fiber above a rational point is commonly connected, in which
case it is the spectrum of a Hurwitz number field.  
 We study many examples
of such maps and their fibers, finding number
fields whose existence contradicts standard mass heuristics.    
\end{abstract}
\maketitle

\tableofcontents

\section{Introduction}   

This paper is a sequel to {\em Hurwitz monodromy and full number fields} \cite{RV15}, joint with Venkatesh.   It is self-contained and aimed more specifically at algebraic number theorists.  
Our central goal is to provide experimental evidence for a conjecture raised in 
\cite{RV15}.  More generally, our objective is to get a concrete and practical
feel for a broad class of remarkable number fields arising in algebraic 
geometry,  the Hurwitz number fields of our title. 

\subsection{Full fields, the mass heuristic, and a conjecture}
\label{ffmh}
    Say that a degree $m$ number field $K=\Q[x]/f(x)$ is {\em full} if the Galois group of $f(x)$ is either the alternating group $A_m$ 
or the symmetric group $S_m$.   For $\cP$ a finite set of primes, let $F_{\cP}(m)$ be 
the number of full fields $K$ of degree $m$ with all primes dividing the discriminant $\disc(K)$ within $\cP$.  

\begin{figure}[htb]
\includegraphics[width=5in]{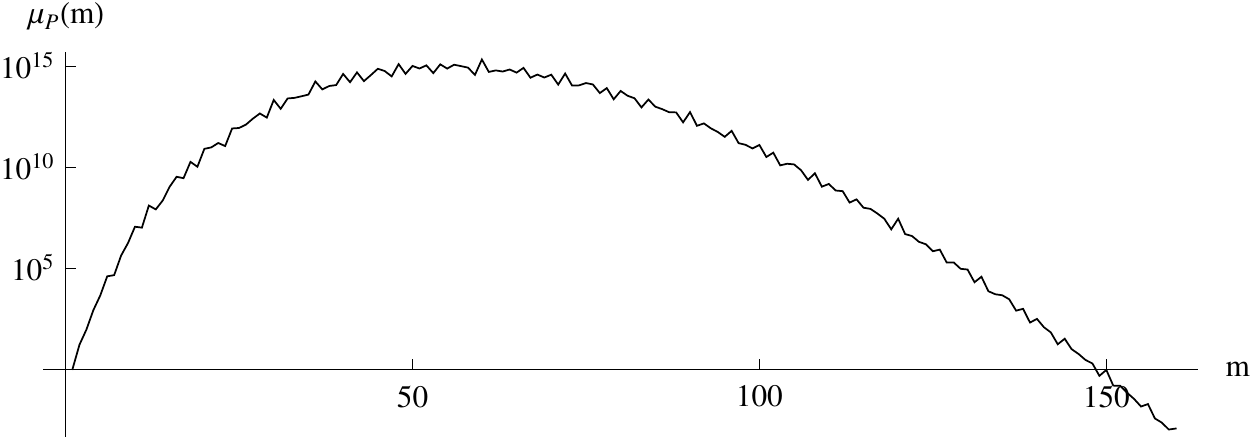}
\caption{\label{mass235}The heuristic approximation $\mu_{\cP}(m)$ to 
the number $F_{\cP}(m)$ of degree $m$ full fields
ramified within $\cP=\{2,3,5\}$}
\end{figure}

In \cite{Bha07}, Bhargava formulated a heuristic expectation $\mu_D(m)$ 
 for the number $F_D(m)$ of degree $m$ full number fields with absolute discriminant $D \in \Z_{\geq 1}$.  The main theorems of \cite{DH71}, \cite{Bha05}, and \cite{Bha10} respectively say that this heuristic 
is asymptotically correct for $m=3$, $4$, and $5$.  While Bhargava is clearly focused in \cite{Bha07}
on this ``horizontal'' direction of fixed $m$ and increasing $D$, it certainly makes
sense to apply the same mass heuristic in the ``vertical'' direction.  
 In \cite{Rob07}, we summed over contributing $D$ to obtain a heuristic expectation $\mu_\cP(m)$
for the number $F_{\cP}(m)$.    Figure~\ref{mass235} graphs the function $\mu_{\cP}$ in
 the completely
typical case of $\cP = \{2,3,5\}$.  For any fixed $\cP$,  the numbers $\mu_\cP(m)$ 
can be initially quite large, but by \cite[Eq.~68]{Rob07}  they ultimately decay super-exponentially
to zero.   From this mass heuristic, one might expect that for 
any fixed $\cP$, the sequence $F_{\cP}(m)$  would be eventually zero.  

The construction studied in \cite{RV15} has origin in work of Hurwitz and
 involves an arbitrary finite
nonabelian simple group $T$.  Let $\cP_T$ be the set of 
primes dividing $T$.   The construction gives a large class of 
separable algebras $K_{h,u}$ over $\Q$ which 
we call {\em Hurwitz number algebras.}  Infinitely many
of these algebras have all their ramification with $\cP_T$.  
Within the range of our computations here, these 
algebras are commonly number fields themselves; in 
all cases, they factor into number fields which we 
call {\em Hurwitz number fields}. The algebras
come in families of arbitrary dimension $\rho \in \Z_{\geq 0}$, 
with the {\em Hurwitz parameter}
$h$ giving the family and the {\em specialization parameter}
$u$ giving the member of the family.     Strengthening Conjecture~8.1
of \cite{RV15} according to the discussion in \S8.5 there, 
we expect that there are enough contributing $K_{h,u}$
to give the following statement. 

\begin{Conjecture} \label{mc} Suppose $\cP$ contains the set of primes
dividing the order of a finite nonabelian simple group. 
Then the sequence $F_{\cP}(m)$ is unbounded.
\end{Conjecture}

\noindent From the point of view of the mass heuristic, the conjecture
has both an unexpected hypothesis and a surprising conclusion.

\subsection{Content of this paper}  

The parameter numbers $\rho=0$ and $1$ have special features 
connected to {\em dessins d'enfants}, and we are presenting families with $\rho \in \{0,1\}$
in \cite{RobHBM}.  To produce enough fields to prove Conjecture~\ref{mc},
it is essential to let $\rho$ tend to infinity.  Accordingly we concentrate here 
on the next case $\rho=2$, with our last example being in the
setting $\rho=3$. 

Section~\ref{deg25a} serves as a quick introduction. 
Without setting up any general framework, it 
exhibits a degree $25$ family.  
Specializing this family gives more than $10,000$ number fields 
with Galois group $S_{25}$ or $A_{25}$ and discriminant of the
form $\pm 2^a 3^b 5^c$.  

Section~\ref{families} introduces Hurwitz parameters
and describes how one passes from a parameter to a Hurwitz
cover.   Full details would require deep forays
into moduli problems on the one hand and braid group
techniques on the other.  We present information 
at a level adequate to provide a framework for 
our examples to come.   In particular, we use the Hurwitz parameter 
$h = (S_5, (2111,5), (4,1))$ corresponding
to our introductory example to illustrate the generalities.

 Section~\ref{specialization} focuses on specialization,
meaning the passage from a Hurwitz cover to its fibers.  
In the alternative language that we have been using in this
introduction,  a Hurwitz cover gives a family
of Hurwitz number algebras, and then specialization 
is passing from the entire family to one of its members.  
The section elaborates on the heuristic argument
for Conjecture~\ref{mc} given in \cite{RV15}.
It formulates Principles
A, B, and C, all of which say that specialization
behaves close to generically.   Proofs of even weak
forms of Principles A and B 
would suffice to prove Conjecture~\ref{mc}.   
Here again, the introductory example 
is used to illustrate the generalities.   

The slightly shorter Sections~\ref{deg9}-\ref{deg1200} 
each report on a family and its specializations,
degrees being $9$, $52$, $60$, $96$, $202$, and $1200$.
Besides describing its family, each section also illustrates 
a general phenomenon.  

Sections~\ref{deg9}-\ref{deg1200} together indicate that the strength with which Principles
A, B, and C hold has a tendency to increase with the degree $m$, in 
strong support of Conjecture~\ref{mc}.   In particular, our two largest degree examples clearly show that Hurwitz 
number fields are not governed by the mass heuristic
as follows.  In the degree $202$ family, Principles A, B, and C 
hold without exception.  
One has $\mu_{\{2,3,5\}}(202) \approx 2 \cdot 10^{-17}$, 
but the family shows 
$F_{\{2,3,5\}}(202) \geq 2947$.
Similarly, $\mu_{\{2,3,5\}}(1200) \approx 10^{-650}$
while the one specialization point
we look at in the degree $1200$ family shows 
$F_{\{2,3,5\}}(1200) \geq 1$.   

There are hundreds of assertions in this paper, with proofs
in most cases involving computer calculations, using 
{\em Mathematica} \cite{math}, {\em Pari} \cite{Pari}, and {\em Magma} \cite{magma}.   We have
aimed to provide an accessible exposition
which should make all the assertions
seem plausible to a casual reader.  We have 
also included enough details so that a diligent
reader could efficiently check any of these
assertions.   Both types of readers could
make use of the large {\em Mathematica} file
\verb@HNF.mma@ on the author's homepage.  
This file contains seven large polynomials 
defining the seven families considered here,
and miscellaneous further information about their
specialization to number fields.  

\subsection{Acknowledgements}  This paper was started at the same time
as \cite{RV15}.  It is a pleasure to thank Akshay Venkatesh whose careful
reading of early drafts of this paper in the context of its relation with \cite{RV15} 
improved it substantially.    It is also a pleasure to thank the Simons foundation
which partially supported this work through grant \#209472.

\section{A degree $25$ introductory family}
\label{sect25a} \label{deg25a}
In this section, we begin by constructing a single  
full Hurwitz number field, of degree $25$ and discriminant $2^{56} 3^{34} 5^{30}$. 
We then use this example to communicate the general nature of Hurwitz number fields and their 
explicit construction.   We close by varying two parameters involved 
in the construction to get more than ten thousand other degree 
twenty-five full Hurwitz number fields from the same family, all ramified within $\{2,3,5\}$.  

\subsection{The  $25$  quintics with critical values $-2$, $0$, $1$ and $2$}
\label{first25}

Consider polynomials in $\C[s]$ of the form
\begin{equation}
\label{fdef}
g(s) = s^5 + b s^3 + c s^2 + d s + x.
\end{equation}
We will determine when the set of critical values of $g(s)$ is $\{-2,0,1,2\}$.  

The critical points of such a polynomial are of course given by the roots 
of its derivative $g'(s)$.   The critical values are then given by the roots of
the resultant
\[
r(t) = \mbox{Res}_s(g(s)-t,g'(s)).
\]
Explicitly, this resultant works out to $r(t) =$
\begin{align*}
&3125 t^4 \\
&+1250 (3 b c-10 x) t^3 \\ 
 &+\left(  108 b^5-900 b^3 d+825 b^2 c^2-11250 b c x+2000 b d^2+2250 c^2 d+18750    x^2   \right) t^2  \\
 &  -2 \left( 108 b^5 x-36 b^4 c d+8 b^3 c^3-900 b^3 d x+825 b^2 c^2  x+280 b^2 c d^2 \right.  \\ 
 & \left. 
    \qquad -315 b c^3 d-5625 b c x^2+2000 b d^2 x+54 c^5+2250 c^2 d x-800 c d^3+6250 x^3
 \right)t \\
  &  + \left(108 b^5 x^2-72 b^4 c d x+16 b^4 d^3+16 b^3 c^3 x-4 b^3 c^2 d^2-900 b^3 d
    x^2+825 b^2 c^2 x^2  \right. \\ 
 & \left. \qquad
    +560 b^2 c d^2 x-128 b^2 d^4-630 b c^3 d x+144 b c^2
    d^3-3750 b c x^3 \right. \\ 
    & \left. \qquad +2000 b d^2 x^2+108 c^5 x-27 c^4 d^2+2250 c^2 d x^2-1600 c
    d^3 x+256 d^5+3125 x^4\right).
\end{align*}
This large expression conforms to the {\em a priori} known structure of $r(t)$: it is a 
quartic polynomial in the variable $t$ depending on the four parameters $b$, 
$c$, $d$, and $x$.    The computation required to obtain the expression is
not at all intensive; for example,  {\em Mathematica}'s \verb@Resultant@
does it in $0.01$ seconds of CPU-time.  

Now consider in general the problem of classifying quintic polynomials \eqref{fdef} with prescribed critical
values.  Clearly, if the given values are the roots of a monic degree four polynomial $\tau(t)$,
then we need to choose the $b$, $c$, $d$, and $x$ so that $r(t)$ is identically equal to 
$3125 \tau(t)$.  
Equating coefficients of $t^i$ for $i=0$, $1$, $2$, and $3$ gives four equations in the four unknowns
$b$, $c$, $d$, and $x$. 
If $(b,c,d,x)$ is a solution then so is $(\omega^2 b, \omega^3 c, \omega^4 d, \omega^5 x)$ for any fifth root of unity $\omega$.    
Thus the solutions come in packets of five, each packet having a common $x$.    

In our explicit example, $\tau(t) = (t+2)t(t-1)(t-2)$.   {\em Mathematica} determines 
in less than a second that there are 125 solutions $(b,c,d,x)$.  
The twenty-five possible $x$'s are the roots of a degree twenty-five polynomial,
\begin{equation}
\label{firstmoduli}
f(x) = 2^{98} 3^8 x^{25} - 2^{96} 3^8 5^2 x^{24} + \cdots + 4543326944239835953052526892234.
\end{equation}
 The algebra $\Q[x]/f(x)$ is our first explicit example
of a Hurwitz number algebra.     In this case, $f(x)$ is irreducible in $\Q[x]$, so 
that $\Q[x]/f(x)$ is in fact a Hurwitz number field.

\subsection{Real and complex pictures}
\label{pictures}
  Before going on to arithmetic concerns, 
we draw two pictures corresponding to the Hurwitz number field $\Q[x]/f(x)$
we have just constructed.  Any Hurwitz number algebra $K$ would have analogous
pictures.    Our object is to visually capture the fact that any 
Hurwitz number algebra $K$ is involved in
a very rich mathematical situation.     Indeed if $K$ has 
degree $m$, then one has $m$ different
geometric objects, with their arithmetic coordinated
by $K$.  

\begin{figure}[htb]
\includegraphics[width=4.6in]{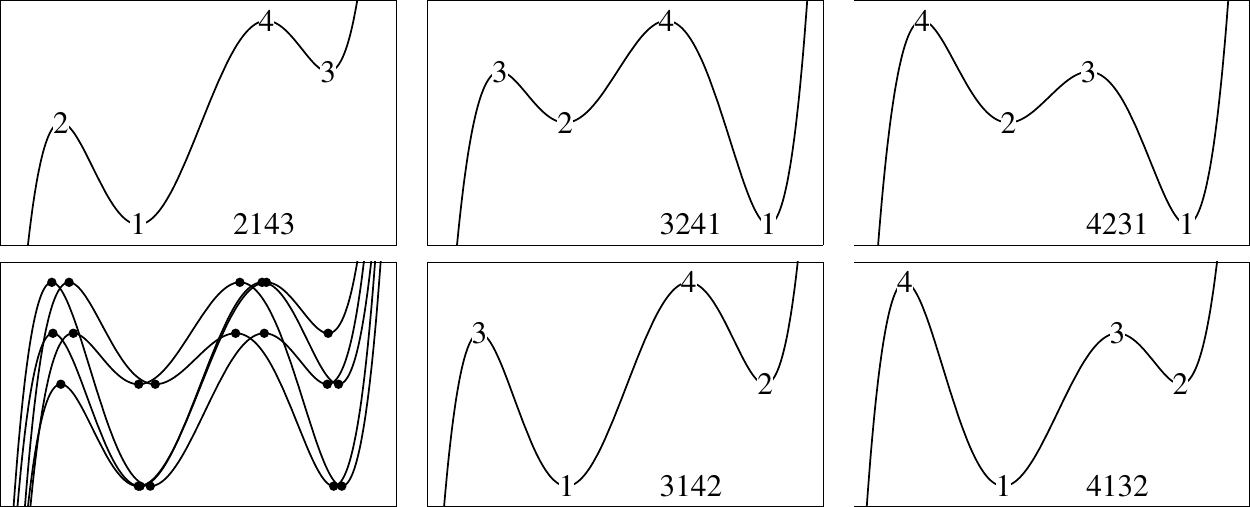}
\caption{\label{realpict} Graphs $t = g_x(s)$ of five quintic polynomials 
with critical points $(s_i,t_i)$ ordered from bottom to top and marked by $i \in \{1,2,3,4\}$. }
\end{figure}
 Of the twenty-five solutions $x$ to \eqref{firstmoduli}, five are real.  Each of these $x$ corresponds to exactly one real solution $(b,c,d,x)$.  The corresponding polynomials
$g_x(s)$ are plotted in the window $[-2.1, 2.1] \times [-2.4,2.4]$ of the real $s$-$t$ plane in  Figure~\ref{realpict}.  
The critical values $t_i$ are indexed from bottom to top so that always $(t_1,t_2,t_3,t_4) = (-2,0,1,2)$, 
with $i$ printed at the corresponding turning point $(s_i,t_i)$.  The labeling of each graph 
encodes the left-to-right ordering of the
critical points $s_i$.  For example, in the upper left rectangle the critical points are 
  $(s_2,s_1,s_4,s_3) \approx (-1.5, -0.6, 0.7, 1.4)$ and the
graph is accordingly labeled by $L=2143$.  The labeling is consistent with 
the labeling in Figure~\ref{pict25} below.     

To get images for all twenty-five roots $x$, we consider the semicircular graph 
$\TopcircleA$
 in the complex $t$-plane
drawn in Figure~\ref{semicircle}.   We then 
 draw in Figure~\ref{pict25} its preimage
  $g_x^{-1}(\TopcircleA)$
  in the complex $s$-plane under
 twenty-five representatives $g_x$.     
 Each of the four critical values $t_i \in \{-2,0,1,2\}$ has a unique
critical preimage $s_i \in \C$, and we print $i$ at $s_i$ in Figure~\ref{pict25}.  
There are braid operations $\sigma_1$, $\sigma_2$, $\sigma_3$ corresponding to 
universal rules which permute the figures, given in this instance by 
Figure~\ref{pict25key}.
Here the $\sigma_i$ all have cycle type $3^5 2^5$ with 
$\sigma_2$ preserving the letter and incrementing the 
index modulo $3$.  The fact that this geometric action has image all
of $S_{25}$ suggests that the Galois group of \eqref{firstmoduli} will 
be $A_{25}$ or $S_{25}$ as well.

{\begin{figure}[htb]
\centering
\begin{minipage}{.05 \textwidth}
$\;$
\end{minipage}
\begin{minipage}{.33\textwidth}
\captionsetup{width=4.0cm}
\includegraphics[width= \textwidth]{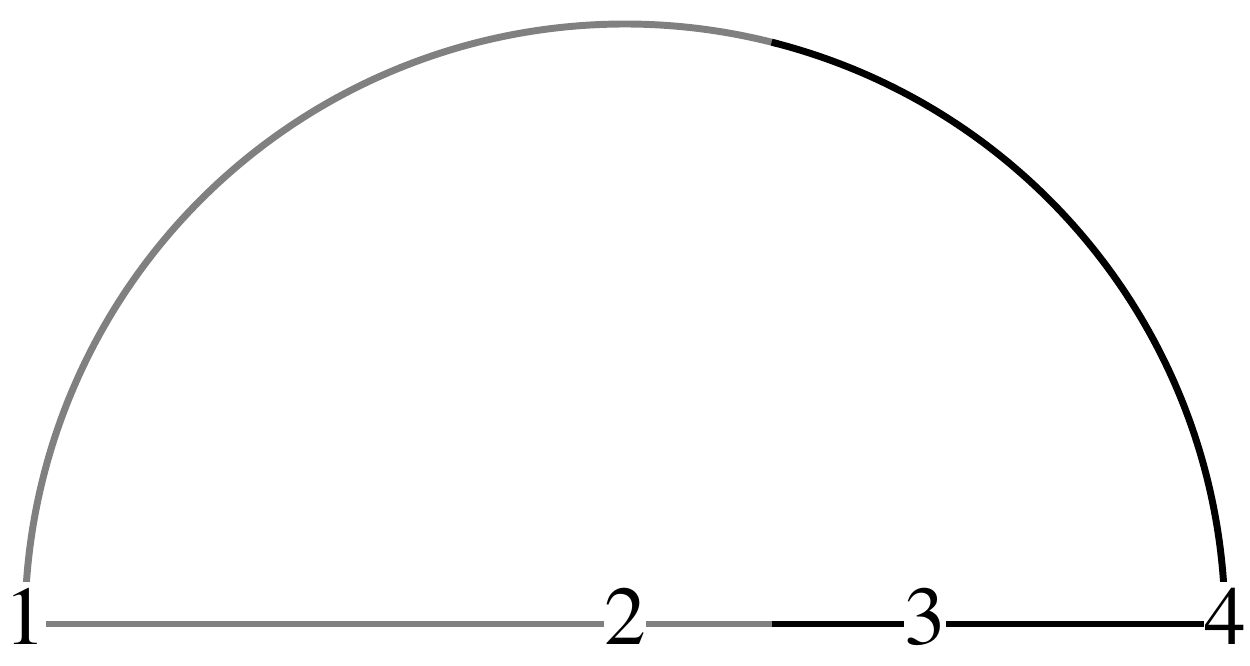}
\caption{\label{semicircle} A 
graph connecting the 
roots of the specialization polynomial 
$(t+2)t(t-1)(t-2)$}
\end{minipage}
\begin{minipage}{.01 \textwidth}
$\;$
\end{minipage}
\begin{minipage}{.57\textwidth}
\captionsetup{width=7.0cm}
 \[\xymatrix@R=.81pc@C=1.2pc{
F \ar[r] \ar[d] & G \ar[r] \ar[d]  & H \ar[d]& A3 \ar[r] \ar[d]  & F \ar[d] \\
B3 \ar[r] \ar[d]  & \ar[r] \ar[d]  I & C3 \ar[d]& B2\ar[r] \ar[d]  & D3 \ar[d] \\
C2 \ar[r] & J \ar[r] & D2 & H \ar[r] & A2 \\
E3  \ar[r] \ar[d] & J \ar[r] \ar[d]  & B1 \ar[d]& E2 \ar[r] \ar[d]  & C1 \ar[d] \\
I \ar[r] & E1 \ar[r] & D1 & G \ar[r] & A1 
}
\] 
\caption{\label{pict25key} Actions of $\sigma_1$ (vertical arrows), $\sigma_2$ (symbols), and $\sigma_3$ (horizontal arrows) }
\end{minipage}
\end{figure}
}

The twenty-five preimages are indeed topologically distinct.    Thus for the twelve  $\gamma_{abc} = \gamma_{cba}$, the critical points $a$, $b$, and $c$ are connected
by a triangle and the middle index $b$ is connected also to the remaining critical point.    
Similarly the indexing for the twelve $\gamma_{abcd} = \gamma_{dcba}$ describes how the 
critical points are connected.   The five graphs corresponding to the real $x$ treated in
Figure~\ref{realpict} are easily identified
by the horizontal line present in Figure~\ref{pict25}. 
\begin{figure}[htb]
\includegraphics[width=4.6in]{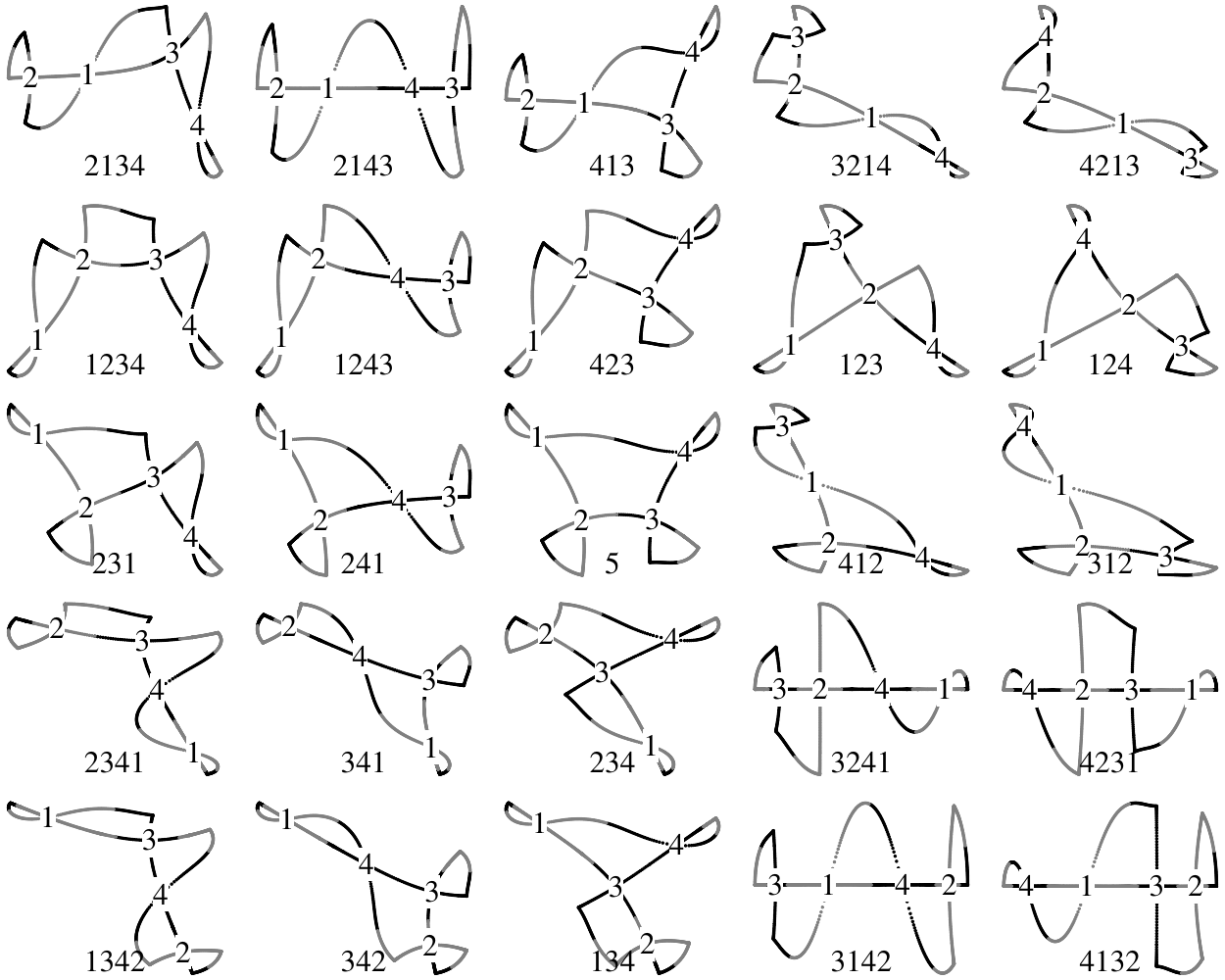}
\caption{\label{pict25} The preimage of the semicircular graph of Figure~\ref{semicircle} under twenty-five quintic polynomials}
\end{figure}  
We are touching on the braid-theoretic infrastructure of Hurwitz number fields in this paper only very lightly.
Our point in presenting Figures~\ref{semicircle}-\ref{pict25} is simply 
to give some idea of the topology behind the existence of Hurwitz number fields.  

\subsection{A better defining polynomial $\phi(x)$ and field invariants.} 
\label{betterpoly}
 We are not so much interested
in the polynomial $f(x)$ from \eqref{firstmoduli} itself, but rather in the field $\Q[x]/f(x)$ it defines.  
{\em Pari}'s command \verb@polredabs@ converts $f(x)$ into a monic polynomial $\phi(x)$ which defines the same field and has minimal sum of the absolute squares of its roots.   It returns
\begin{eqnarray*}
\lefteqn{\phi(x)=} \\
&&  x^{25}-5 x^{24}+15 x^{23}-5 x^{22}-380 x^{21}+1290 x^{20}-4500 x^{19}-28080
    x^{18} \\ && +183510 x^{17}+74910 x^{16}-3033150 x^{15}+4181370 x^{14}+27399420
    x^{13}\\&& -48219480 x^{12}-124127340 x^{11}+266321580 x^{10}+466602765
    x^9\\&&  -592235505 x^8  -905951965 x^7+1232529455 x^6+2423285640 x^5\\&& +664599470
    x^4-814165000 x^3 -517891860 x^2-58209720 x+2436924.
\end{eqnarray*}
For fields of sufficiently small degree, one applies the reduction operation \verb@polredabs@
 as a matter of course: 
the new smaller-height polynomials are more reflective of the complexity of the fields considered, isomorphic fields may be revealed, and any subsequent analysis of field invariants is sped up.

{\em Pari}'s \verb@nfdisc@ calculates that the discriminant of $\Q[x]/\phi(x)$ is
\[
D=1119186718586212624367616000000000000000000000000000000=2^{56} 3^{34} 5^{30}.
\]
The fact that $D$ factors into the form $2^a 3^b 5^c$ is known from the 
general theory presented in Sections~\ref{families} and \ref{specialization}, using that $2$, $3$, and $5$ are the primes less than or equal to the degree $5$, and the polynomial discriminant of $\tau(t)=(t+2)t(t-1)(t-2)$, namely $2304=2^8 3^2$, has this form too. 
 Note that since all the 
exponents of the field discriminant are greater than the degree $25$, the 
number field is wildly ramified at all the
base primes, $2$, $3$, and $5$. 

 To look more closely at $\Q[x]/\phi(x)$, 
we factorize the $p$-adic completion
$ \Q_p[x]/\phi(x)$
as a product of fields over $\Q_p$.  We write the symbol $e^f_c$ to indicate a factor
of degree $ef$, ramification index $e$, and discriminant $p^{fc}$.   One gets
\begin{align}
\nonumber \mbox{$2$-adically:} \;\;\;\;\;  & {\bf 16}_{50} \, 3_2 \,  3_2 \, 3_2,  \\
\label{fieldsymbols} \mbox{$3$-adically:} \;\;\;\;\;  & {\bf 9}_{18} \, 4_3 \, 4_3  \, {\bf 3}_5 \, {\bf 3}_5 \, 1^2_0, \\
\nonumber \mbox{$5$-adically:} \;\;\;\;\;  & {\bf 25}_{30},
\end{align}
with wild factors printed in bold.   Thus, the first line means that
$\Q_2[x]/\phi(x)$ is a product $K_1 \times K_2 \times K_3 \times K_4$,
where $K_1$ is a wild totally ramified degree sixteen extension of $\Q_2$ with discriminant $2^{50}$, while $K_2$, $K_3$,  and $K_4$
are tame cubic extensions of discriminant $2^2$. 
The behavior for the three primes is roughly typical, 
although, as we'll see in Figure~\ref{discfchart}, a little less ramified than average.

Because the field discriminant is a square, the  Galois group of $\phi(x)$ is in $A_{25}$.  Many small collections of $p$-adic factorization patterns for small unramified $p$ each suffice to 
prove that the Galois group 
is indeed all of $A_{25}$.   Most easily, $\phi(x)$ factors in $\Q_{19}[x]$ into irreducible factors
of degrees $17$, $6$, and $2$, so that the Galois group contains an element of
order $17$.  Jordan's criterion now applies: a transitive subgroup of $S_m$ containing  \index{criterion, Jordan's}
an element of prime order in $(m/2,m-2]$ is all of $A_m$ or $S_m$.   We will use this easy technique without further comment for
 all of our other determinations that
 Galois groups of number fields are full.   
 One could also use information from ramified primes as above, 
 but unramified primes give the easiest computational route.

\subsection{A family of degree $25$ number fields} 
\label{fam25two}
We may ask, more generally, for the quintics with any fixed set of critical values.
  This amounts to repeating our previous computation, replacing the polynomial
  $\tau(t) = (t+2)t(t-1)(t-2)$ of the three previous subsections 
with other separable quartic polynomials 
\begin{equation}
\label{original}
\tau(t) =  t^4 + b_1 t^3 + b_2 t^2 + b_3 t + b_4.
\end{equation}
From each such $\tau$, we obtain a degree $25$ algebra over $\Q$, once again
the algebra determined by the possible values of the variable $x$.  

Changing $\tau$ via a rational affine transformation $t \rightarrow \alpha t + \beta$ does not  change the degree twenty-five algebra constructed.  
Accordingly, one can restrict attention to 
specialization polynomials $\tau(t)$ with $b_1=0$, 
and consider only a set of representatives for the equivalence 
$(b_2, b_3, b_4) \sim (\alpha^2 b_2,  \alpha^3 b_3, \alpha^4 b_4)$, where
$\alpha$ is allowed to be in $\Q^\times$.  
 In particular, if $b_2$ and $b_3$ are nonzero, any such polynomial
 is equivalent to a unique polynomial   of the form 
\begin{equation}
\label{normalized}
\tau(u,v,t) = t^4-2 t^2 v-8 t v^2-4 u v^2+v^2.
 \end{equation}
 Here the reason for the complicated form on the right is explained in the discussion around \eqref{discagree}. 
 We will treat in what follows only the main two-parameter family where $b_2$ and $b_3$ are both nonzero.
Note, however, that two secondary one-parameter families are also interesting:
if $b_3=0$, one gets degree $25$ algebras with
Galois group in 
$S_5 \times S_2 \wr S_{10}$, because of the symmetry induced from
$t \mapsto -t$;  the case $b_2=0$ gives rise to 
full degree $25$ algebras, just like the main case.  

One can repeat the computation of \S\ref{first25}, now with the parameters $u$ and
$v$ left free.   The corresponding general degree twenty-five moduli polynomial $f_{25}(u,v,x)$ has $129$ terms as an expanded polynomial in $\Z[u,v,x]$.  
After replacing $x$ by $5x/4$ and clearing a constant, coefficients average about $16$ digits. 
We will not write this large polynomial explicitly here, instead giving a simpler polynomial that 
applies only in the special case $u=1/3$ at the end of \S\ref{realpictures}.
 
 \subsection{Keeping ramification within $\{2,3,5\}$} 
 \label{keeping} Suppose $\tau(t)$ from 
  \eqref{original} normalizes to $\tau(u,v,t)$ from \eqref{normalized}.    We write the corresponding
  Hurwitz number algebra as $K_{u,v}$.   Inclusion \eqref{Pinclude} below says that if $\tau(t)$ is ramified within $\cP  = \{2,3,5\}$,
  then so is $K_{u,v}$.   By a computer search we have found $11031$ 
  such $(u,v)$.    From irreducible $f_{25}(u,v,x)$, we obtain
  $F_\cP(25) \geq 10938$.  The remaining $f_{25}(u,v,x)$ all
  have a single rational root and from these polynomials we obtain
  $F_\cP(24) \geq 93$.   The behavior of the $11031$ different $K_{u,v}$ will 
  be discussed in more detail 
  in Section~\ref{specialization} below.
  
  A point to note is that ramification is obscured by the passage to standardized
  coordinates.  In the case of our first example $\tau(t) = (t+2)t(t-1)(t-2)$, 
  the corresponding $(u,v)$ is $(37/175,9/1715)$.   The standardized polynomial
   $\tau(37/175,9/1715,t)$ after 
  clearing denominators has a $7$ in its discriminant. 
  

\section{Background on Hurwitz covers}
 \label{families}
In this section, we provide general background on Hurwitz covers.
  Most of our presentation is in the setting
of algebraic geometry over the complex numbers.  In the 
last subsection, we shift to the more arithmetic setting
where Hurwitz number fields arise.  

\subsection{Hurwitz parameters}  
\label{parameters}
We use the definition in \cite[\S1B]{RV15} of Hurwitz parameter:
  {\em Let $r \in \Z_{\geq 3}$. 
 An $r$-point {\em Hurwitz parameter} is a triple $h = (G,C,\nu)$ 
where 

$\bullet \;$ $G$ is a finite group;

$\bullet \;$ $C = (C_1,\dots,C_k)$ is a list of conjugacy classes whose union generates $G$;

$\bullet \;$ $\nu = (\nu_1,\dots,\nu_k)$ is a list of positive integers summing to $r$ such that 
$\prod [C_i]^{\nu_i} = 1$ in the abelianization $G^{\rm ab}$.}
%
 We henceforth always take the $C_i$ distinct and not the identity, and
 normalize so that $\nu_i \geq \nu_{i+1}$.
  The number $\nu_i$ functions
as a multiplicity for the class $C_i$.  

Table~\ref{7families} gives the Hurwitz parameters of the seven 
Hurwitz covers described in this paper.  It is also gives the associated degrees
$m$ and bad reduction sets $\cP_h$, each to be discussed later in this section.
\begin{table}[htb]
\[
\begin{array}{r| lll | rl}
\mbox{Section} & G & C & \nu & m & \cP_h \\
\hline
\mbox{\S\ref{deg25a}}
 & S_5 & (2111, 5) & (4,1) &  25 & \{2,3,5\}   \\
\mbox{\S\ref{deg9}} & S_3 \wr S_2 & (21111, 33, 222) & (3,1,1) & 9 & \{2,3\}  \\
\mbox{\S\ref{deg52}} & S_6 & (21111,222,3_1111,3_\infty2_01) & (2,1,1,1) & 52 & \{2,3,5\} \\
\mbox{\S\ref{deg60}} & PSL_3(3) & (2^4 1^5,3^3 1^4) & (3,2) & 2 \cdot 60 & \{2,3\} \\
\mbox{\S\ref{deg96}} & GL_3(2) & (22111, 421) & (4,1) & 2 \cdot 96 & \{2,3,7\} \\
\mbox{\S\ref{deg202}} & S_6 & (21111, 3_021,3_1111,4_\infty11) & (2,1,1,1) & 202 & \{2,3,5\} \\
\mbox{\S\ref{deg1200}} & S_6 & (21111, 321, 411) & (4,1,1) & 1200 &  \{2,3,5\}
\end{array}
\]
\caption{\label{7families} Hurwitz parameters for the seven covers pursued in this paper, two of them 
with normalizations given via subscripts}
\end{table}
In the case that $G$ is a symmetric group $S_n$,  
we label a conjugacy class $C_i$ by the partition $\lambda_i$ of $n$ giving
the lengths of the  cycles of any of its elements.    We describe classes for general 
$G$ in a similar way.   Namely we choose a transitive embedding
 $G \subseteq S_n$.  We then label classes $C_i$ by their induced cycle
 partitions $\lambda_i$, removing any ambiguities which arise by further labeling.  
 In none of our examples is further labeling necessary.   

 Our concept of Hurwitz parameter emphasizes multiplicities more than other similar
concepts in the literature.  For example, the first line of Table~\ref{7families} says that 
our introductory example comes from the parameter
$h = (S_5, (2111, 5),(4,1))$.    In e.g.\ \cite{MM99}, the indexing scheme
would center on the class vector $(2111, 2111, 2111, 2111,5)$.   

\subsection{Covers indexed by a parameter} 
\label{covers} An $r$-point parameter $h = (G,C,\nu)$ determines
an unramified cover of $r$-dimensional complex algebraic varieties 
\begin{equation}
\label{cover1}
\pi_h : \AHur_h \rightarrow \AConf_\nu.
\end{equation}
 The base
 is the variety whose points are tuples $(D_1,\dots,D_k)$ of disjoint subsets $D_i$ 
 of  the complex projective line $\AP^1$, with $D_i$ consisting of $\nu_i$ 
 points.   Above a point $u = (D_1,\dots,D_k) \in \AConf_\nu$, the fiber has one 
 point for each solution of a moduli problem indexed by $(h,u)$.   
   
 The  moduli problem described in \cite[\S2]{RV15}
 involves degree $|G|$ Galois covers $\Sigma \rightarrow \AP^1_t$, with 
 Galois group identified with $G$.   An equivalent
  version of this moduli problem 
  makes reference to the embedding
 $G \subseteq S_n$ used to label conjugacy classes.  
 When $G$ is its own normalizer in $S_n$, 
 which is the case for all our examples,
 the equivalent version is
 easy to formulate:
 above a point $u = (D_1,\dots,D_k)  \in \AConf_\nu$, the fiber $\pi_h^{-1}(u)$
 consists of points $x$ indexing isomorphism classes of degree $n$ covers 
 \begin{equation}
 \label{YtoP}
\AS_x \rightarrow \AP_t^1.  
 \end{equation}
 These covers are required to have global monodromy group $G$, local 
 monodromy class $C_i$ for all $t \in D_i$, and be otherwise unramified.    
 In this equivalent version, the ramification numbers of the preimages
 of $t \in D_i$ in $\AS_x$ together form the partition $\lambda_i$.  
 
We prefer the equivalent version for the purposes of this 
paper, since it directly guides our actual computations.  
For example, in our introductory example, the
quintic polynomials prominent there can
 be understood as degree five rational 
 maps $\AP^1_s \rightarrow \AP^1_t$.  Here 
 $\AP^1_s$ is a common coordinatized version of all the  $\AS_x$.  
 Also the preimage of $\infty$ consists 
 of the single point $\infty$, explaining 
 why polynomials rather than more
 general rational functions are involved.
 At no point did degree $120$ maps explicitly enter into the computations of  
 Section~\ref{deg25a}.
 
 \subsection{Covering genus}
Let $h = (G,C,\nu)$ be a Hurwitz parameter with $G \subseteq S_n$ a transitive permutation group.
Let $\ell_i$ be the the number of parts of the partition $\lambda_i$ induced by $C_i$, and let $d_i = n-\ell_i$ be the corresponding drop. 
Consider the Hurwitz covers $\AS_x \rightarrow \AP^1_t$ parametrized by $x \in \AHur_h$.   By the Riemann-Hurwitz formula,
the curves $\AS_x$ all have genus $g = 1 -n +  \frac{1}{2} \sum \nu_i d_i$.   

Given $G$, let $d$ be the minimal drop of a nonidentity element.  If $h$ is an $r$-point Hurwitz
parameter based on $G$, then necessarily $g \geq 1-n+dr/2$.  To support Conjecture~\ref{mc},
one needs to draw fields from cases with arbitrarily large $r$ and thus 
arbitrarily large $g$.   However explicit computation of families
rapidly becomes harder as $g$ increases, and in this paper we only
pursue cases with genus zero.

 \subsection{Normalization}  
 \label{normalization}
 The three-dimensional complex group $\APGL_2$ acts by fractional
 linear transformations on $\AConf_\nu$.  Since $\APGL_2$ is connected, the
 action lifts uniquely to an action on $\AHur_h$ making $\pi_h$
 equivariant.    To avoid redundancy, it is important for us to use this action to replace
 \eqref{cover1} by a cover of varieties of dimension $\rho = r-3$.   Rather than working
 with quotients in an abstract sense, we work with explicit codimension-three
 slices as follows.  
  
We say that a Hurwitz parameter is {\em base normalizable} if 
 $k \geq 3$ and $\nu_{k-2}=\nu_{k-1}=\nu_k=1$.  For
 a base normalizable Hurwitz parameter, we replace 
  \eqref{cover1} by a map of
 $\rho$-dimensional varieties,
 \begin{equation}
\label{cover1a}
\pi_h : \AX_h \rightarrow \AU_\nu. 
 \end{equation}
 Here the target $\AU_\nu$ is the subvariety of 
 $\AConf_\nu$ with $(D_{k-2},D_{k-1},D_k) = (\{0\},\{1\},\{\infty\})$.  
 The domain $\AX_h$ is just the preimage of $\AU_\nu$ in 
 $\AHur_h$.   This reduction in dimension is ideal for our purposes:
 each $\APGL_2$ orbit on $\AConf_\nu$ contains 
 exactly one point in $\AU_\nu$.  
 
 We say that a base normalizable genus zero Hurwitz parameter
 is {\em fully normalizable} if the partitions $\lambda_{k-2}$, $\lambda_{k-1}$, and $\lambda_{k}$ have between them 
 at least three singletons.  A normalization is then obtained  
%
 %
%
 by labeling
 three of the singletons by $0$, $1$, and $\infty$, as illustrated
 twice in Table~\ref{7families}.   This labeling places a unique coordinate function $s$ 
 on each $\AS_x$.  Accordingly, each point of $\AX_h$ is then identified 
 with an explicit rational map from $\AP_s^1 \rightarrow \AP_t^1$.  
 

 When the above normalization conventions do not apply, we modify 
 the procedure, typically in a very slight way, so as to likewise replace the cover of 
 $r$-dimensional varieties \eqref{cover1} by a cover
 of $\rho$-dimensional varieties \eqref{cover1a}.  For example,
 two other multiplicity vectors $\nu$ figuring into some of our examples
are $(4,1)$ and $(3,1,1)$.   For these cases, we define
\begin{align*}
\tau_4(t) & = t^4-2 t^2 v-8 t v^2-4 u v^2+v^2, &
\tau_3(t) & = t^3 + t^2 + u t + v.
\end{align*} 
The form for $\tau_4(t)$ is chosen to make
discriminants tightly related:
\begin{align}
\label{discagree}
\disc_t(\tau_4(t)) & =  -2^{12} v^6 d, &
\disc_t(t \tau_3(t)) & =   v d, 
\end{align}
with
\begin{equation}
\label{quartdisc}
d=4 u^3-u^2-18 u v+27 v^2+4 v.
\end{equation}
In the respective cases, we say that a divisor tuple is normalized
if it has the form 
\begin{align*}
(D_1,D_2) & = ((\tau_4(t)),\{\infty\}), &
(D_1,D_2,D_3) & = ((\tau_3(t)), \{0\}, \{\infty\}).
\end{align*}
These normalization conventions define subvarieties 
$\AU_{4,1} \subset \AConf_{4,1}$ and 
$\AU_{3,1,1} = \AConf_{3,1,1}$.  
As explained in the $(4,1)$ setting in \S\ref{fam25two}, we are throwing away 
some perfectly interesting $\APGL_2$ orbits on $\AConf_\nu$
by our somewhat arbitrary normalization conventions.  However all these
orbits together have positive codimension in $\AConf_\nu$ and what 
is left is adequate for our purposes of supporting Conjecture~\ref{mc}.  
Always, once we have $\AU_\nu \subset \AConf_\nu$ we just take $\AX_h \subset \AHur_h$ to be
its preimage.

The two base varieties just described are identified by their common
coordinates: $\AU_{4,1} = \AU_{3,1,1} = \Spec \,\C[u,v,1/vd]$.   
This exceptional identification has a conceptual source as
follows.  With $(u,v)$ fixed, let $D_1 = (\tau_4(u,v,t))$ so that $(D_1,\{\infty\}) \in \AU_{4,1}$.    Let $V$ be the four-element 
subgroup of $\APGL_2$ consisting of fractional transformations
stabilizing the roots of $\tau_4(u,v,t)$.  One then has a degree four map $q$  from $\AP^1_t$ to
its quotient $\AP :=\AP^1_t/V$.   There are three natural divisors
on $\AP$:  the divisor $\Delta_1$ consisting of the three critical values,
and the one-point divisors $\Delta_2 = \{q(D_1)\}$ and $\Delta_3 = \{q(\infty)\}$. 
Uniquely coordinatize $\AP$ so that 
$(\Delta_1,\Delta_2,\Delta_3) = ((\tau_3(u',v',t)),\{0\},\{\infty\})$.  Then $u'=u$ and
$v'=v$.

 

\subsection{The mass formula and braid representations}

The degree $m$ of a cover $\AX_h \rightarrow \AU_\nu$ can 
be calculated by 
%
 group-theoretic techniques
as follows.  Define the {\em mass} $\overline{m}$ of an $r$-point Hurwitz parameter $h = (G,C,\nu)$
via a sum over the irreducible characters of $G$:
\begin{equation}
\label{massformula}
\overline{m}=  \frac{\prod_i |C_i|^{\nu_i}}{|G|^2}
 \sum_{\chi \in \widehat{G}}  \frac{\prod_i \chi(C_i)^{\nu_i}}{\chi(1)^{r-2}}.
\end{equation}
Then $\overline{m} \geq m$ always.  
Suppose that no proper subgroup $H \subset G$ contains elements from
all the conjugacy classes $C_i$,
as is the case in \S\S\ref{deg25a}, \ref{deg9}, \ref{deg52}.  
Then $\overline{m}=m$.  When there
 are exist such $H$, as in \S\S\ref{deg60}, \ref{deg96}, \ref{deg202}, and \ref{deg1200},
 one can still get exact degrees by applying \eqref{massformula}
 to all such $H$ and computing via inclusion-exclusion.  Chapter~7 of \cite{Se08} gives \eqref{massformula}
 and works out several examples in the setting $r=3$. 
 
As a one-parameter collection of examples, 
consider $h(j) = (S_5, (2111,5), (j,1))$ for $j \geq 4$ even.  
Since $S_5$ is generated by any $5$-cycle and 
any transposition, one has $\overline{m} = m$ for $h(j)$.  
From $0$'s in the character table of $S_5$, only the characters $1$, $\epsilon$, $\chi$, and $\chi \epsilon$ contribute, 
with $\epsilon$ the sign character and $\chi+1$ the given degree $5$ permutation character.   
 We can ignore $\epsilon$ and
$\chi \epsilon$ by doubling the contribution of $1$ and $\chi$:
\begin{eqnarray*}
m = \frac{10^j  24}{120^2} \left(2  + 2   \frac{\chi(2111)^j \chi(5)}{\chi(1)^{r-1}} \right) = 
 \frac{10^{j-2}}{6} \left(2 + 2 \cdot  \frac{2^j (-1)}{4^{j-1}}  \right) = \frac{1}{3} \left( 10^{j-2} - 5^{j-2} \right).
\end{eqnarray*}
For $j=4$, one indeed has $m=25$, as in the introductory example.

The monodromy group of a cover $\AX_h \rightarrow \AU_\nu$ 
can be calculated by group-theoretic techniques \cite[\S3]{RV15}.  These techniques
center on braid groups and underlie the mass formula.   The output of these
calculations is a collection of permutations in $S_m$ which generate the 
monodromy group, with $\langle \sigma_1,\sigma_2,\sigma_3\rangle = S_{25}$ from Figure~\ref{pict25key}
being completely typical.   Fullness 
of these representations is important for us: once we switch over
to the arithmetic setting in \S\ref{rationality},  it implies fullness
of generic specializations. 

Theorem~5.1 of \cite{RV15} proves
a general if-and-only-if result about fullness.  In one
direction, the important fact for us here  is that  to systematically
obtain fullness one needs for $G$ to be very close 
to a nonabelian simple group $T$.  Here ``very close'' includes subgroups
of $\Aut(T)$ of the form $T.2$, such as $G=S_n$ for $T = A_n$.  
 This direction
accounts for the hypothesis of Conjecture~\ref{mc}.
In the other direction, fullness is the 
typical behavior for these $G$.  This statement is the main theoretical
reason we expect that the conclusion of 
Conjecture~\ref{mc} follows from the hypothesis.

%



\subsection{Accessible families}
\label{accessible1} 
 The groups $A_n$ and $S_n$ give rise to many computationally accessible 
 families with $\rho \in \{0,1,2\}$.
 Table~\ref{accessible} presents families
with $\rho=2$ and $n \in \{5,6\}$, omitting $1$'s from partitions to save space.
The table gives the complete list 
of $h$ with covering genus $g=0$ and degree $m \in \{1,\dots,250\}$.  We 
have verified by a braid group computation
that the $58$ families listed
all have full monodromy group.     
\begin{table}[htb]
{\small
{\renewcommand{\arraycolsep}{1.5pt}
\[
\begin{array}{| r | lllll | r|}
                                       \hline
                     & \multicolumn{5}{c|}{\mbox{$\nu$ = (3,1,1)}}&\\
                     n &  C_1 & C_1 & C_1 & C_2 & C_3 & m  \\
                     \hline
                     6&3  &3  &3  &2  &32 &216 \\
                     6&2  &2  &2  &32 &5  &150 \\
                     6&2  &2  &2  &32 &42 &120 \\
                     6&3  &3  &3  &2  &4  &96  \\
                     6&2  &2  &2  &4  &5  &75  \\
                     6&2  &2  &2  &4  &42 &72  \\
                     6&22 &22 &22 &2  &222&60  \\
                     6&2  &2  &2  &22 &6  &C54  \\
                     6&2  &2  &2  &32 &33 &C54  \\
                     5&2  &2  &2  &3  &4  &48  \\
                     5&2  &2  &2  &22 &4  &48  \\
                     5&2  &2  &2  &3  &32 &45  \\
                     6&3  &3  &3  &2  &222&44  \\
                     6&2  &2  &2  &3  &6  &36  \\
                     6&2  &2  &2  &4  &33 &B36  \\
                     5&2  &2  &2  &22 &32 &B36  \\
                     \bf{6}&\mathbf{2}  &\mathbf{2}  &\mathbf{2}  &\mathbf{222} &\mathbf{5}  &\mathbf{A25}  \\
                                           \hline
                                                  \multicolumn{7}{c}{\;} \\
       \multicolumn{7}{c}{\;} \\
       \multicolumn{7}{c}{\;} \\
                        \end{array}
\;\;\;\;\;\;\;\; 
\begin{array}{| r | lllll | r |}
                     \hline
                     & \multicolumn{5}{c|}{\mbox{$\nu$ = (2,1,1,1)}}&\\
                     n &  C_1 & C_1 & C_2 & C_3 & C_4 & m  \\
                     \hline
                     6 &  22  & 22  & 2   & 3   & 4   & 240    \\
                    {\bf 6} &  {\bf 2}   & {\bf 2}   & {\bf 3}  &{\bf  4}   & {\bf 32}  & {\bf 202}   \\
                     6 &  3   & 3   & 2   & 22  & 4   & 168    \\
                     6 &  2   & 2   & 3   & 22  & 5   & 125    \\
                     6 &  2   & 2   & 3   & 22  & 42  & 100   \\
                     6 &  2   & 2   & 22  & 32  & 222 & 60     \\
                     6 &  22  & 22  & 2   & 3   & 222 & 57     \\
                     {\bf 6} &  {\bf 2}   & {\bf 2}   & {\bf 3}   & {\bf 32}  & {\bf 222} & {\bf 52}     \\
                     6 &  2   & 2   & 3   & 22  & 33  & 48    \\
                     6 & 3   & 3   & 2   & 22  & 222 & 42     \\
                     6 &  2   & 2   & 22  & 4   & 222 & 40     \\
                     6 &  2   & 2   & 3   & 4   & 222 & 36    \\
                     \hline
                          \multicolumn{7}{c}{\;} \\
                      \hline
                     & \multicolumn{5}{c|}{\mbox{$\nu$= (4,1)}}&\\
                     n &  C_1 & C_1 & C_1 & C_1 & C_2 & m  \\
                     \hline
                     6&3  &3  &3  &3  &22 &192 \\
                     {\bf 5}&{\bf 2}  &{\bf 2}  &{\bf 2}  &{\bf 2}  &{\bf 5}  &{\bf A25}  \\
                       \hline
           \multicolumn{7}{c}{\;} \\
       \multicolumn{7}{c}{\;} \\
       \multicolumn{7}{c}{\;} \\
\end{array}
\;\;\;\;\;\;\;\; 
\begin{array}{| r | lllll | r |}
                      \hline
                     & \multicolumn{5}{c|}{\mbox{$\nu$= (2,2,1)}}&\\
                     n &  C_1 & C_1 & C_1 & C_2 & C_3 & m \\
                     \hline
                     6 &  2   & 2   & 22  & 22  & 5   & 175   \\
                     6 &  2   & 2   & 4   & 4   & 22  & 158   \\
                     6 &  2   & 2   & 22  & 22  & 42  & 128   \\
                     6 &  2   & 2   & 4   & 4   & 3   & 89     \\
                     6 &  2   & 2   & 3   & 3   & 42  & 80     \\
                     6 &  2   & 2   & 3   & 3   & 5   & 75    \\
                     6 &  2   & 2   & 22  & 22  & 33  & 54  \\
                     5 &  2   & 2   & 3   & 3   & 22  & 58    \\
                     5 & 2   & 2   & 22  & 22  & 3   & 48     \\
                     6 &  2   & 2   & 3   & 3   & 33  & 39    \\
                     \hline
                                                                                    \multicolumn{7}{c}{\;} \\
                      \hline
                     & \multicolumn{5}{c|}{\mbox{$\nu$ = (5)}}&\\
                     n & C_1 & C_1 & C_1 & C_1 & C_1 & m  \\
                     \hline
                     6&3  &3  &3  &3  &3  &96  \\
                                    \hline
                                                              \multicolumn{7}{c}{\;} \\
                      \hline
                     & \multicolumn{5}{c|}{\nu=(3,2)}&\\
                     n &  C_1 & C_1 & C_1 & C_2 & C_2 & m  \\
                     \hline
                            5&3  &3  &3  &2  &2  &55  \\
                            5&22 &22 &22 &2  &2  &40  \\
                       \hline
                        \multicolumn{7}{c}{\;} \\                     
\end{array}
\]
}
}
\caption{\label{accessible} Fifty-eight computationally accessible two-parameter families.
One, eight, one, and forty-eight of these families respectively have $G = A_5$, $S_5$, $A_6$, and $S_6$.
 }
\end{table}

Table~\ref{accessible} reveals that our introductory example has the lowest degree
$m$ in this context.  It and the only other degree $25$ family are highlighted in bold.  
Two of the six families we pursue in \S\ref{deg9}-\ref{deg1200} are likewise put 
in bold.  The remaining families from these sections are not on the table because three of them
have group different from $A_n$ and $S_n$ and one has $\rho=3$.  

A remarkable phenomenon revealed by braid computations is what we call {\em cross-parameter agreement.}    There are three instances on Table~\ref{accessible}: covers given with the same
label, be it $A$, $B$, or $C$, are isomorphic.   Note that the first instance involves
the exceptional isomorphism $\AU_{4,1} = \AU_{3,1,1}$ from \S\ref{normalization}, with the cover of 
$\AU_{4,1}$ being our introductory family.   Many instances of cross-parameter agreement are
given with defining polynomials in \cite{RobHBM}.   V\"olklein \cite{Vol01} explains some instances of 
cross-parameter agreement via the Katz middle convolution operator \cite{Kat96}.    

\subsection{Computation and rational presentation} 
\label{computation}
Our general method of passing from a Hurwitz parameter $h = (G,C,\nu)$ to an explicit Hurwitz cover is well
illustrated by our introductory example.  Very briefly, one writes down all covers $\AS \rightarrow \AP^1_t$
conforming to $h$ and satisfying the chosen normalization conditions.  From this first step,
one extracts a generator $x$ of the function field of the  variety $\AX_h$.  
For all $\nu$ we are considering, one has also coordinates
 $u_1$, \dots, $u_\rho$ on the base variety $\AU_\nu$.
By computing critical values, one arrives at a degree $m$ polynomial relation $f(u_1,\dots,u_\rho,x)=0$
describing the degree $m$ extension $\C(\AX_h)/\C(\AU_\nu)$.   In all the 
examples of both Table~\ref{accessible} and \S\ref{deg9}-\ref{deg1200}, the covering variety $\AX_h$ is connected 
and so $\C(\AX_h)$ is a field.  In general, as illustrated many times in 
\cite{RobHBM}, the polynomial $f(u_1,\dots,u_\rho,x)$ may factor, making $\AX_h$ disconnected and $\C(\AX_h)$ a product of fields.  

When $\AX_h$ is a connected rational variety, one can seek a more insightful presentation
as follows.  One finds not just the above single element $x$ of the function field,
but rather elements $x_1$, \dots, $x_\rho$ which satisfy $\C(\scX_h) = \C(x_1,\dots,x_\rho)$.
Then, working birationally, the map $\pi_h : \AX_h \rightarrow \AU_\nu$ is given
by $\rho$ rational functions,
\begin{equation}
\label{rationalpresentation}
u_i = \pi_{h,i}(x_1,\dots,x_\rho).
\end{equation}
We call such a system a {\em rational presentation}.

As an example of a rational presentation, consider the Hurwitz parameter $\hat{h}_{25} = (S_6,(21111,222,51),(3,1,1))$,
chosen because it relates to our introductory example $h_{25}$ by cross-parameter
agreement.  We partially normalize via
$5_\infty 1_0$.  
We complete
our normalization by requiring the coefficient of $s^2$ in the cubic in the numerator of $g(s)$ be $1$: 
\begin{align*}
g(s) & = \frac{\left(s^3+s^2+ z s + y\right)^2}{as}, &\frac{g'(s)}{g(s)}  & = \frac{5 s^3+3 s^2+z s - y}{s \left(x+s^3+s^2+s z\right)}.
\end{align*}
In the logarithmic derivative of $g(s)$ to the right, let $\Delta(s)$ be its numerator.  
Writing $g(s) = g_0(s)/g_\infty(s)$, one requires that the resultant $\mbox{Res}_s(g_0(s)-g_\infty(s) t,\Delta(s))$ be
proportional to $t^3 + t^2 + u t + v$.  Working out this proportionality
makes $a =  4(27 - 225z + 500z^2 + 375y - 5625yz)/3125$.   

We have thus identified $\AX_h$ birationally with the plane $\C_y \times \C_z$.  But moreover, the proportionality gives
\begin{eqnarray}
\label{uh25} u & = & \frac{5^5 \left(\!\! \!\!\!\!\!\!\!\! \begin{array}{c}-2025 y^3+2700 y^2 z^2-405 y^2 z-12 y^2-660 y z^3 \\ \qquad  \qquad \qquad \qquad 
+301 y z^2-36 y z+16 z^5-8
    z^4+z^3 \end{array} \!\!  \right)}{\left(-5625 y z+375 y+500 z^2-225 z+27\right)^2}, \\
\label{vh25}  v & = & -\frac{5^{10} y \left(27 y^2-18 y z+4 y+4
    z^3-z^2\right)^2}{\left(-5625 y z+375 y+500 z^2-225 z+27\right)^3}.
 \end{eqnarray}
 Equations~\ref{uh25} and \ref{vh25} together form a rational presentation of the form \eqref{rationalpresentation}.  
 In general, one can always remove all but one of the $x_i$ by resultants, thereby returning to a $\rho$-parameter
 univariate polynomial.

    To see the cross-parameter agreement between $h_{25}$ and $\hat{h}_{25}$ explicitly,
    we proceed as in \cite[(5.3) or (5.5)]{RobG2} to identify the root of $f_{25}(u,v,x)$ in the function
    field $\C(y,z)$.  It turns out to be 
   \begin{equation}
   \label{agree}
   x =  \frac{3 \cdot 5^7 z \left(4 y^3-y^2-18 y z+27 z^2+4 z\right)}{2 \left(500
    y^2-5625 y z-225 y+375 z+27\right)^2}.
 \end{equation}
  Thus the natural function $x$ in the first approach has
  only a rather complicated presentation in the second approach.  

 \subsection{Rationality, descent, and bad reduction}  
 \label{rationality}
 We have been working over $\C$ so far in this section to 
 emphasize that large parts of our subject matter are 
 a mixture of complex geometry and group theory.  
 In the construction of Hurwitz number fields, arithmetic 
 enters ``for free'' and only at the end.   For example,
 the final equations  \eqref{uh25} and
 \eqref{vh25} have coefficients in $\Q$, even though
 we were thinking only in terms of complex 
 varieties when deriving them.  
 


Following \cite[\S2D]{RV15} we say that a Hurwitz 
parameter $h = (G,C,\nu)$ is {\em strongly rational}
if all the conjugacy classes $C_i$ are rational.  This is the case in all our examples,
as each $C_i$ is distinguished from 
all the other classes in $G$ by its partition $\lambda_i$.
We henceforth work only with strongly rational Hurwitz parameters.    
In this case,
the cover \eqref{cover1} canonically descends
to a cover of varieties defined over $\Q$,
\begin{equation}
\label{cover2}
\pi_h : \Hur_h \rightarrow \Conf_\nu.
\end{equation}
Similarly, since all our normalizations are chosen rationally, 
the corresponding reduced cover \eqref{cover1a} descends to a cover of $\Q$-varieties,
$\pi_h : \scX_h \rightarrow \scU_\nu$. Computations 
as in our introductory example or the previous subsection
end at polynomials $f(u_1,\dots,u_\rho, x) \in \Q[u_1,\dots,u_\rho,x]$
whose vanishing corresponds to \eqref{cover2}.

 Note that in the previous paragraph we changed 
fonts as we passed from complex spaces to $\Q$-varieties.   
As a further example of this font change, $\AU_\nu$ has appeared many 
times already as conveniently brief notation 
for $\scU_\nu(\C)$.   In the future we will also need the
subsets  $\scU_\nu(R)$ for 
various subrings $R$ of $\C$.    In subsequent sections we
will continue this convention: when working primarily geometrically
we emphasize complex spaces, and when specializing 
we emphasize varieties over $\Q$.


Let $\cP_h$ be the set of primes at which
\eqref{cover2} has bad reduction.  Let $\cP_G$ 
be the set of primes dividing the order of $G$.
Then a fundamental fact is
\begin{equation}
\label{Pinclude}
\cP_h \subseteq \cP_G.  
\end{equation}
This fact is essential for our argument supporting Conjecture~\ref{mc},
and enters our considerations through \eqref{keyproperty}.    
The inclusion \eqref{Pinclude} follows from the standard reference \cite{BR11} because
all the results there hold for any ground field with
characteristic not dividing $|G|$.  
Table~\ref{7families} gives $\cP_h$ for our covers.  

\section{Specialization to Hurwitz number algebras}
\label{specialization}  This section discusses 
specializing a given Hurwitz cover $\scX_h \rightarrow \scU_\nu$
to number fields, taking the introductory example of Section~\ref{deg25a} further
to illustrate general concepts.   The goal is to extrapolate from the
observed behavior of the $11031$ algebras $K_{u,v} = K_{(S_5,(2111,5),(4,1)),(u,v)}$ to
the expected behavior of specialization in general.  We dedicate a subsection 
each to Principles A, B, and C.  The extent to which they hold 
will be discussed in connection with all of our examples in
the sequel.  

\subsection{Algebras corresponding to fibers}
Let $\scX_h \rightarrow \scU_\nu$ 
be a Hurwitz cover, as in \S\ref{rationality}.   Let $u \in \scU(\Q)$.
The scheme-theoretic fiber $\pi_h^{-1}(u)$
is the spectrum of a $\Q$-algebra $K_{h,u}$.  We call 
$K_{h,u}$ a Hurwitz number algebra.  The homomorphisms of
$K_{h,u}$ into $\C$ are indexed by points of the complex fiber 
$\pi_h^{-1}(u) \subset \AX_h$.
Like all separable algebras,
the $K_{h,u}$ are products of fields.  These factor fields are the Hurwitz
number fields of our title.  Whenever the monodromy group of
$\AX_h \rightarrow \AU_\nu$ is transitive, the algebras $K_{h,u}$ 
are themselves fields for generic $u$, by the Hilbert irreducibility 
theorem.   

For many $\nu$, certainly including all $\nu$ containing three $1$'s,
 $\scU_\nu$ can be identified with an open subvariety of affine space
$\Spec \,\Q[u_1,\dots,u_\rho]$ as in \cite[\S8]{Rob15}.  Birationally at least, the cover is given by a polynomial
equation $f(u_1,\dots,u_\rho,x) = 0$.   The point $u$ corresponds
to a vector $(u_1,\dots,u_\rho) \in \Q^\rho$.  The algebra $K_{h,u}$
is then $\Q[x]/f(u_1,\dots,u_\rho,x)$.   The factorization of
$K_{h,u}$ into fields corresponds to the factorization of $f(u_1,\dots,u_\rho,x)$ 
into algebras.

\subsection{Real pictures and specialization sets $\scU_\nu(\Z[1/\cP])$} 
\label{realpictures}
\begin{figure}[htb]
\includegraphics[width=4.8in]{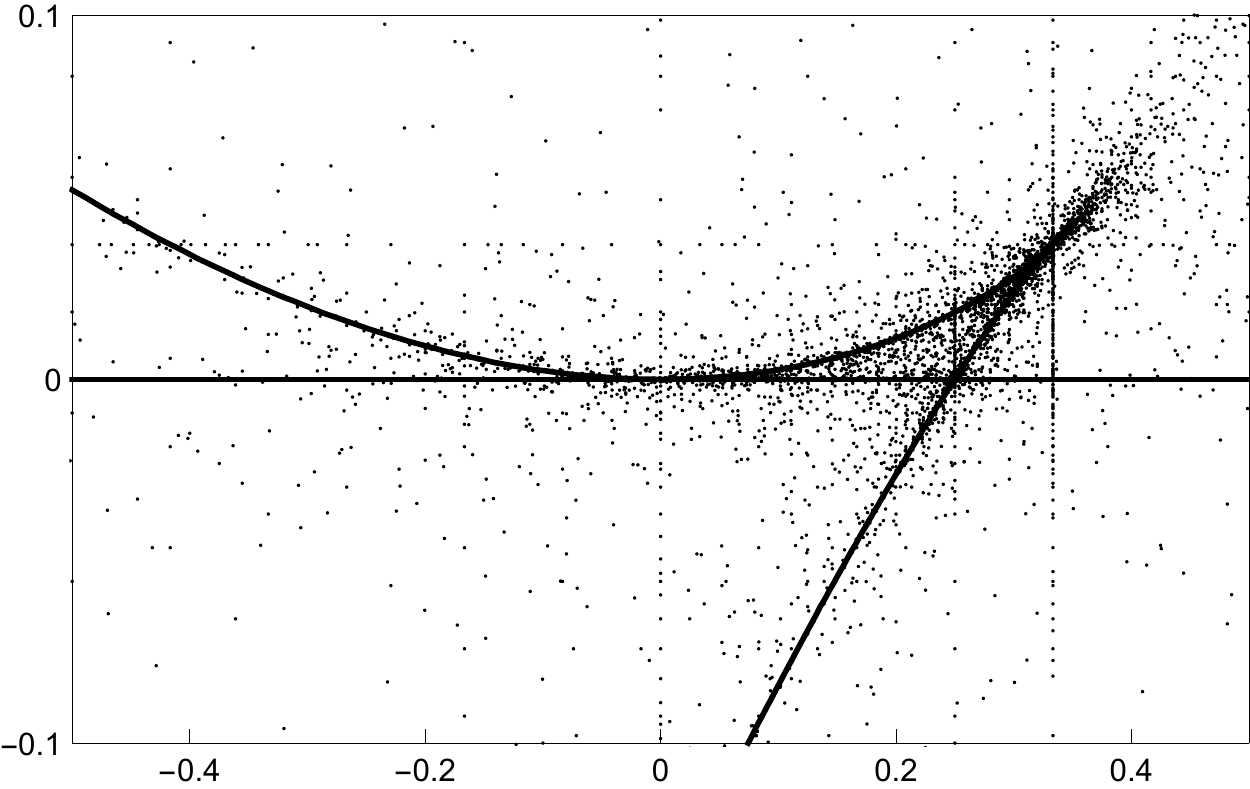}
\caption{\label{scatterintro}  A window on $U_{4,1}(\R)$, which is the complement 
of the two discriminantal curves in the $u$-$v$ plane.  Points are part of the
specialization set $\scU_{4,1}(\Z[1/30])$, which applies in \S\ref{deg25a}.   
For \S\ref{deg9}, \S\ref{deg60}, and \S\ref{deg96}, specialization
sets $\scU_{3,1,1}(\Z[1/6])$, $\scU_{3,1,1}(\Z[1/6])$, and $\scU_{4,1}(\Z[1/42])$
are respectively used, and the corresponding pictures would
have the same discriminant locus but different specialization points.  }
\end{figure}

Figure~\ref{scatterintro} draws a window on $\scU_{4,1}(\R)$. 
With the choice of coordinates made in \S\ref{normalization}, it is the complement
of the drawn discriminant locus in the real $u$-$v$ plane. One
should think of the line at infinity in the projectivized plane
 as also part of the discriminant locus.
An analogous picture for $\nu=(2,1,1,1)$ is
drawn in Figure~\ref{spec2111}.

Let $\cP$ be a finite set of primes with product $N$.  Let $\Z[1/\cP] = \Z[1/N]$
be the ring obtained from $\Z$ by inverting the primes in $\cP$.   When 
the last three entries of $\nu$ are all $1$ then $\scU_\nu$ is naturally
a scheme over $\Z$.   Accordingly it makes sense
to consider $\scU_\nu(R)$ for any commutative ring. 
 The finite set of points $\scU_\nu(\Z[1/\cP])$ is
studied in detail in \cite{Rob15}, including complete identifications
for many $(\nu,\cP)$.   For general $\nu$, one similarly has a finite
subset $\scU_\nu(\Z[1/\cP])$ of $\scU_\nu(\Q)$.    Its key property for us is that 
\begin{equation}
\label{keyproperty}
\parbox{4in}{{\em for any Hurwitz cover $\scX_h \rightarrow \scU_\nu$ and 
any $u \in \scU_\nu(\Z[1/\cP])$, the algebra $K_{h,u}$ is 
ramified within $\cP_h \cup \cP$.}}
\end{equation}
In \S\ref{ramificationtame} we take $\cP$ strictly containing
$\cP_h$ so as to provide examples of ramification known {\em a priori} to be 
tame.  Otherwise, we are always taking $\cP = \cP_h$ in this paper.  
Figure~\ref{scatterintro} shows the 8461 of the known $11031$ points of  $\scU_{4,1}(\Z[1/30])$ which fit into the 
window.

In our Hurwitz parameter formalism, we are emphasizing the multiplicity vector
$\nu$ because of the following important point.  
Fix $r$ and a non-empty  finite set of primes $\cP$, and 
consider all multiplicity vectors $\nu$ with total $r$.  Then 
$\scU_\nu(\Z[1/\cP])$ tends to get larger as 
$\nu$ moves from $(1^r)$ to $(r)$.   
This phenomenon is represented
by the two cases considered for $\cP = \{2,3,5\}$ in
this paper: $|\scU_{2,1,1,1}(\Z[1/30])| = 2947$, from \cite[\S8.5]{Rob15},
and $|\scU_{4,1}(\Z[1/30])| \geq 11031$.  
In fact, as $r$ increases the cardinality $|\scU_{1^r}(\Z[1/\cP])|$ 
eventually becomes zero \cite[\S2.4]{Rob15} while 
$|\scU_{r-3,1,1,1}(\Z[1/\cP])|$ increases without bound \cite[\S7]{Rob15}.
This increase is critical in supporting Conjecture~\ref{mc}.       

In both Figure~\ref{scatterintro} and the similar Figure~\ref{spec2111}, one can see specialization
points from $\scU_\nu(\Z[1/30])$ concentrating on certain lines.  These lines, and other less visible
curves, have the property that they intersect the discriminant locus in the projective plane
exactly three times.  
While the polynomial $f_{25}(u,v,x)$ of \S\ref{fam25two} was too complicated to print, variants
over any of these curves are much simpler.  For example, the most prominent of the
lines is $u=1/3$.  Parametrizing this line by $v = (j-1)/27j$, one has the simple equation
\begin{eqnarray*}
\lefteqn{
f_{25}(j,x)  =  2^2 (x+2) \cdot } \\ && \left(729 x^8-486 x^7-702 x^6-8 x^5+105
    x^4+1118 x^3-1557 x^2+1296 x-576\right)^3 \\
   && \qquad  + 5^{15}  j (x-1)^4 x^9.
\end{eqnarray*}
The ramification partitions above $0$, $1$, and $\infty$ are respectively $3^81$, $2^{10} 1^5$, 
and $(12,9,4)$.   A systematic treatment of these special curves in the cases $\nu=(3,1,1)$
and $\nu = (3,2)$ is given in \cite[\S7]{RobG2}.  For general $\nu$, 
they play an important role in
\cite{RobHBM}.  In this paper the above line $v=1/3$ will play a prominent role
 in \S\ref{deg96}, and analogous lines for $\nu = (2,1,1,1)$ will enter in \S\ref{specsixcurves} and \S\ref{view}.

\subsection{Pairwise distinctness}  
 For each of the $11031$ algebras $K_{u,v}$ of \S\ref{keeping}, and each prime $p \geq 7$,
 one has a Frobenius partition $\alpha_{u,v,p}$ giving the degrees of the factor fields of 
$K_{u,v} \otimes \Q_p$.   For $p=7$, $11$, $13$, $17$, $19$, and $23$, the number
of partitions of $25$ arising is $71$, $126$, $157$, $205$, $243$, and $302$. 
Taking now  $p=7$, $11$, $13$, $17$, $19$, and $23$ as 
cutoffs, the number of tuples $(\alpha_{u,v,7},\dots,\alpha_{u,v,p})$
arising is $71$, $2992$, $10252$, $10981$, $11027$, and $11031$. 
Thus the $11031$ algebras are pairwise non-isomorphic.  
There are many other quick ways of seeing this pairwise 
distinctness.  For example, 
one could use that $6772$ different discriminants 
$D_{u,v}$  arise as a starting point.  

Abstracting this simple observation to a general Hurwitz map $\scX_h \rightarrow \scU_\nu$ gives
\medskip

\noindent{\bf Principle A.}  {\em For almost all pairs of distinct elements $u_1$, $u_2$ in
$\scU_\nu(\Z[1/\cP])$, the algebras $K_{h,u_1}$ and $K_{h,u_2}$ 
are non-isomorphic. }
\medskip 

\noindent So, at least when one restricts to the known $11031$ elements of $\scU_{4,1}(\Z[1/30])$, 
Principle~A holds without exception for our introductory family.  

\subsection{Minimal Galois group drop }  
The Galois group of $f_{25}(u,v,x)$ over $\Q(u,v)$ is $S_{25}$.  
Some of the $11031$ specialized algebras $K_{u,v}$ have smaller 
Galois groups as follows.   First, in $93$ cases, there is a factorization of the form 
$K_{u,v} = K_{u,v}' \times \Q$, with $K_{u,v}'$ a field.   Second, the discriminant
of the specializing polynomial $\tau(u,v, t)$ and the discriminant of the degree twenty-five algebra $K_{u,v}$ agree
modulo squares.  Thus one knows the total number of 
times that a given discriminant class
$d \in \Q^\times/\Q^{\times 2}$ occurs, 
even without inspecting the $K_{u,v}$ themselves.
The number of degree $m$ fields obtained with discriminant class
$d$ is as follows:
\[  
{\renewcommand{\arraycolsep}{2.4pt}
\begin{array}{r|rrrrrrrrrrrrrrrr}
   &   -30 &-15&-10&-6 &-5 &-3  &-2  &-1 &1  &2  &3  &5  &6   &10 &15 &30\\
\hline
25 & 1050 & 547 & 310 & 363 & 641 & 1702 & 1000 & 480 & 557 & 360 & 576 & 572 &
      1026 & 787 & 897 & 70 \\
24  &  14 & 3 & 2 & 4 & 5 & 15 & 8 & 4 & 2 & 4 & 10 & 6 & 3 & 1 & 12 & 0
\end{array}.
}
\]
Galois groups are as large as possible given the above considerations.  Thus
  $A_{25}$ and $A_{24}$ occur respectively $557$ times and twice, 
 leaving $S_{25}$ and $S_{24}$ occurring respectively 
 $10381$ and $91$ times.
 
 To state a principle for general $h = (G,C,\nu)$, let $\Gal_h$ 
 be the generic Galois group of the cover $\scX_h \rightarrow \scU_\nu$.
 
 \medskip

\noindent{\bf Principle B.}  {\em For almost all elements $u$ in
$\scU_\nu(\Z[1/\cP])$, the specialized Galois group $\Gal(K_{h,u})$
contains the derived group $\Gal'_h$ of the generic Galois group.  }
\medskip 

\noindent The most important case of this principle for us is 
when $\scX_h \rightarrow \scU_\nu$ is full, i.e.\ all of $A_m$ or $S_m$.
Then the principle says that $K_{h,u}$ is full for almost
all $u \in \scU_\nu(\Z[1/\cP])$.   In our example, $93$ of the $11031$
known points of $\scU_{4,1}(\Z[1/30])$, thus slightly less than $1\%$, are exceptions to 
the principle.  However, these exceptions
are relatively minor, in that they produce contributors
to $F_{\{2,3,5\}}(24)$ rather than $F_{\{2,3,5\}}(25)$.  

Principle B is formulated so that it includes other cases 
of interest to Conjecture~\ref{mc}.   For example, 
let $m = m_1+m_2$ with $m_1,m_2 \geq 3$.    Suppose 
$\Gal_h$ is one of the five intransitive groups 
containing $A_{m_1} \times A_{m_2}$.  Then Principle B holds
for $u$ if and only if $K_{h,u}$ factors as a product
of two full fields.  This case is illustrated many times
in \cite{RobHBM}, with splittings of the form $25=10+15$ and $70=30+40$
being presented in detail in \S7.1 and \S7.2 respectively.

\subsection{Wild ramification.}  
\label{wildram} 
\begin{figure}[htb]
\includegraphics[width=4.6in]{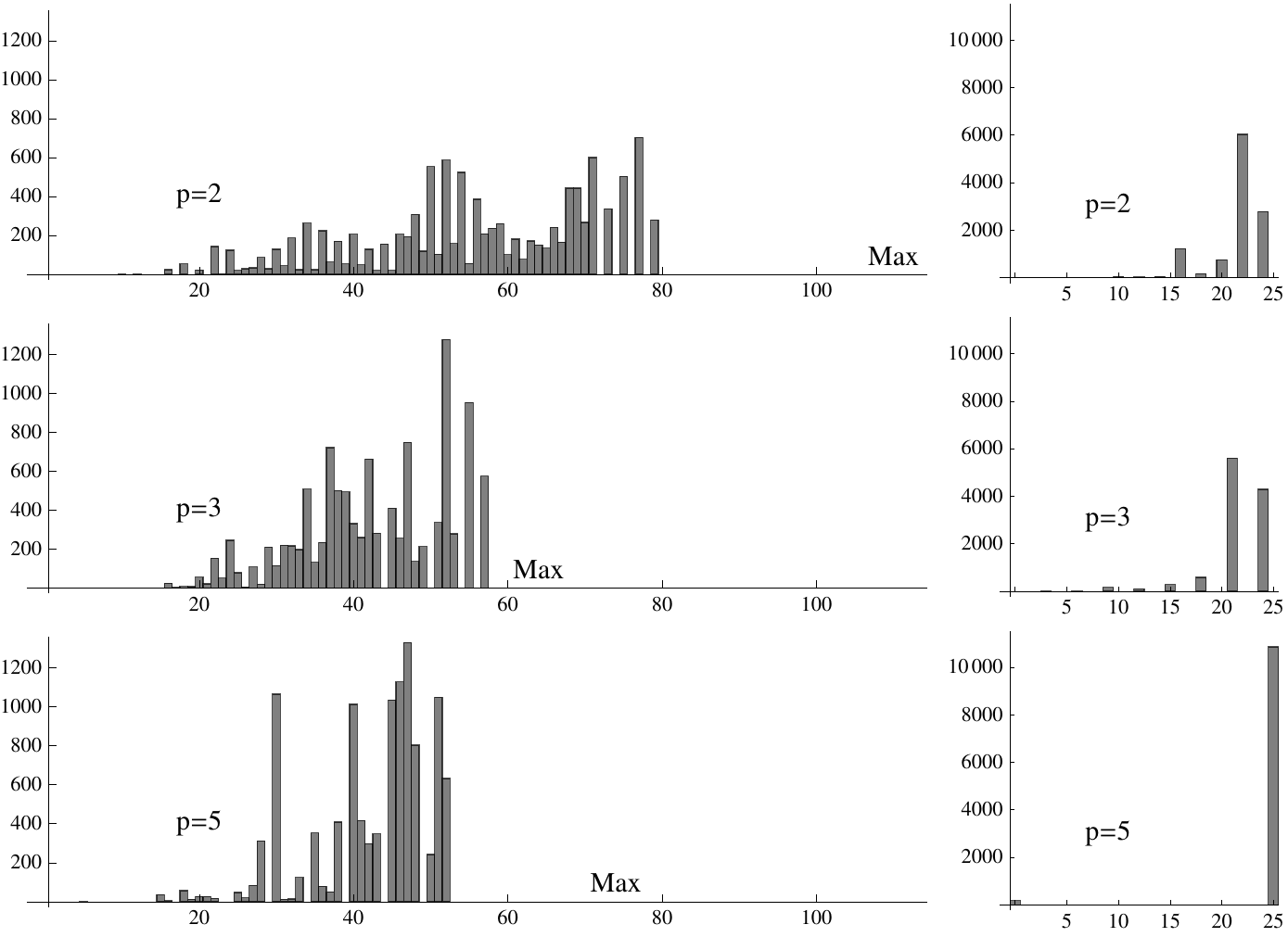}
\caption{\label{discfchart}  Left: distribution of the discriminant exponents  $\ord_p(D)$ the 11031 algebras $K_{u,v}$; the variation of $\ord_p(D)$ is much less than is allowed by general discriminant bounds.  Right: distribution of the wildness degrees
 $m_{p{\rm\mbox{-}wild}}$ 
relevant for Principle C.}
\end{figure}

Consider the discriminants
$\disc(K_{u,v}) = \pm 2^a 3^b 5^c$ as $(u,v)$ varies over the 11031 known elements of $\scU_{4,1}(\Z[1/30])$.
The left part of 
 Figure~\ref{discfchart} gives  the distribution of the exponents $a$, $b$ and 
 $c$.  There is much less variation in the exponents than is allowed for 
field discriminants of degree twenty-five algebras in general.  For general algebras, 
the minimum value 
for $a$, $b$, and $c$ is of course $0$ in each case.  The maximum values occur for the algebras defined
by 
$(x^{16}-2)(x^8-2)x$, $(x^{18}-3)(x^6-3)x$, and $x^{25}-5$, and are respectively $110$, $64$,
and $74$.    The average values in our 
family are $(\langle a \rangle, \langle b \rangle, \langle c \rangle) \approx (56,43,42)$.  

There are many open questions to pursue with regard to wild ramification.  
One could ask for lower bounds valid for all $u$, upper bounds valid for all $u$, or even exact formulas for wild ramification as a function of $u$.    
 Principle C is 
in the spirit of lower bounds:
\medskip 

\noindent{\bf Principle C.} {\em For almost all $u \in \scU_\nu(\Z[1/\cP])$, the specialized algebra
$K_{h,u}$ is wildly ramified at all primes $p \in \cP_h$.}
\medskip

\noindent An interesting complement to Principle C is a general upper bound: $K_{h,u}$ can be
wildly ramified at $p$ only if $p \in \cP_h$ or $p \leq \max(\nu_i)$.  


Certainly if $\ord_p(K_{h,u}) \geq m$ then Principle C holds for $K_{h,u}$ and $p$.  
The left part of Figure~\ref{discfchart} shows that, for each 
$p$, most $K_{u,v}$ satisfy this sufficient criterion.  In fact, for $p=2$, $3$, and $5$,
there are only 374, 568, and 179 algebras $K_{u,v}$ which do not.   However 
to conform to Principle C at $p$, an algebra $K_{h,u}$ needs only to satisfy a much
weaker condition.   Define the {\em wild degree} of a $p$-adic algebra $K$ to
be the sum of the degrees of its wildly ramified factor fields.  Thus in \eqref{fieldsymbols} 
these degrees $m_{p\rm{\mbox{-}wild}}$ for $p=2$, $3$, and $5$ are $16$, $15$, and
$25$ respectively.   Then conformity to Principle~C at $p$ means simply that the $p$-adic 
wild degree is positive.    

The right part of  Figure~\ref{discfchart} gives the distribution of $m_{p\rm{\mbox{-}wild}}$. 
In the cases $p=2$, $3$, and $5$, the heights of the bars above multiples of $p$ are
 \begin{align*}
 & (0, 0, 2, 0, 0, 45, 20, 44, 1211, 161, 746, 6029, 2773), \\
 & (0, 4, 12, 164, 95, 284, 585, 5598, 4289),  \\
 & (179, 0, 0, 0, 0, 10852). 
 \end{align*}
 Thus all 11031 algebras are wild at $2$ and $3$, most having
 $m_{p\rm{\mbox{-}wild}}$ near its maximum possible value of $24$.   
 For $p=5$, there are $179$ exceptions to Principle C, including all the 93 factorizing
 algebras. Besides these exceptions, all algebras have  $m_{5\rm{\mbox{-}wild}}$
 at its maximum possible value of $25$.  
 
 \subsection{Expectations}  As discussed in \cite[\S8]{RV15}, the Hilbert irreducibility theorem already points
 in the direction of  Principles A and B.   In a wide variety of contexts, analogs of these principles hold with great strength.
 For example, in \cite[\S9]{MR05} several covers are discussed in the setting $\cP = \{2,3\}$ and for 
 most of them both Principles A and B hold without exception.   However the situation we are 
 considering here, with fixed $\cP$ and arbitrarily large degree $m$, is outside the realm of previous
 experience.   Explicitly verifying the principles in degrees large enough to contradict
 the mass heuristic
  is important for being confident that these
 standard expectations do indeed hold in this new realm.
 
We are confident that for a given $G$ and varying $h = (G,C,\nu)$, one has 
 strict inclusion $\cP_h \subset \cP_G$ for only finitely many $(C,\nu)$.  
 This expectation, together with Principle C, suggests that there are only
 finitely many full fields $K_{h,u}$ ramified strictly within $\cP_G$.   
One possibility is that full number fields coming from Hurwitz-like constructions are 
the main source of outliers to the mass heuristic.   If one believes this, then
one is led to the first of the two extreme possible complements to
 Conjecture~\ref{mc} discussed
at the end of \cite{RV15}:
%
%
 {\em The sequence $F_{\cP}(m)$ 
 always has support on a density zero set, and it is eventually zero
 unless $\cP$ contains the set of primes divisors of the order of
 a nonabelian finite simple group}.    Our verification that Principle
 C holds with great strength in our examples is supportive
 of this very speculative assertion.  
 

 
 



 \section{A degree 9 family: 
 comparison with complete number field tables}
\label{fam9}  \label{deg9}  
 This section begins our sequence of six sample families of increasing
 degree.  To start in very low degree, we take $G$ solvable. 
 The number fields coming from this first example
are not full and so not directly relevant to
Conjecture~\ref{mc}. This family is nonetheless a good place
to begin our presentation of examples, as
the low degree makes 
 comparison with complete tables of number fields possible.  
 
 \subsection{A Hurwitz parameter with solvable $G$}   Consider
 first the Hurwitz parameter $h' = (S_6, (21111,222,33),(3,1,1))$.
 The mass formula \eqref{massformula} yields $\overline{m}_{h'} = 9$.
 In the language of \cite[\S3]{RV15}, 
 $S_6$ has nine orbits on $\cG_{h'}$.    However
 none of the tuples $(g_1,\dots,g_5) \in \cG_{h'}$ 
 generate $S_6$. 
 
 To place this degenerate situation into our formalism, let 
 $G$ be the wreath product $S_3 \wr S_2$ of order $72$,
 considered as a subgroup of $S_6$.  
 The group $G$ has unique conjugacy classes with cycle type $21111$, $222$, and $33$.
 In place of $h'$, take $h = (G,(21111,222,33),(3,1,1))$.   This 
 parameter accounts for everything, as 
 $\overline{m}_h = m_h=9$.  
 

\subsection{A two-parameter polynomial.} In the present context of $\nu = (3,1,1)$, our normalized specialization polynomials take the form 
\[
\tau(u,v,t) = (t^3 + t^2 + u t + v) t.
\]
The discriminant of the cubic factor is
$d=
 4 u^3-u^2-18 u v+27 v^2+4 v$ from \eqref{quartdisc}.     
A nonic polynomial capturing the family and a resolvent octic 
are as follows: 
\begin{eqnarray*}
f_9(u,v,x) & = & 
x^9-3 x^8+12 u x^7-4  (u+12 v) x^6 +42 v x^5   -6 (4 u+1) v x^4\\ 
&& \qquad +4 v
     (2 u+3 v) x^3 -12 v^2 x^2+3 (4 u-1) v^2 x-v^2 (4 u-8 v-1), \\
&&\\
f_8(u,v,x) & = & x^8+x^4 \left(18 v-6 u^2\right)+x^2 \left(8 u^3-36 u v+108
    v^2\right) \\ 
    && \qquad +(-3 u^4+18 u^2 v-27 v^2).
\end{eqnarray*}
Here $f_9(u,v,x)$ and $f_8(u,v,x)$ respectively have Galois group $9T26 = \F_3^2.GL_2(\F_3)$ and 
$8T23 = GL_2(\F_3) = \tilde{S}_4$.  Because of the complete lack of singletons in the partitions 222 and 33, our 
computation of $f_9(u,v,x)$ required substantial {\em ad hoc} deviations from the procedure
sketched in \S\ref{computation}.

The discriminants of the two polynomials are respectively
\begin{align*}
D_9(u,v) & =  -2^{24} 3^9 v^{10} d^4  (27 v-1)^6, &
D_8(u,v) & =  -2^{24} 3^{19} v^8 d^4 \left(u^2-3 v\right)^2.
\end{align*}
In each case, the discriminant modulo squares is $-3$.  Because of this constancy, 
the Galois groups of $f_9(u,v,x)$ and $f_8(u,v,x)$  over $\C(u,v)$ are respectively
the index two subgroups  
$9T23=\F_3^2.SL_2(\F_3)$ and $8T12 = SL_2(\F_3) = \tilde{A}_4$.  
The last factor of the discriminant
 in each case is an artifact of our particular
polynomials; these factors do not contribute to field discriminants in specializations.

\subsection{Comparison of specializations with complete tables of number fields}
We work with $507$ pairs $(u,v)$ in $\SU_{3,1,1}(\Z[1/6])$.    Twenty-one of them have
$v=1/27$ and so $f_9(u,v,x)$ is not separable.   For forty more, $f_9(u,v,x)$ 
also reduces, with the factorization partitions $63$, $81$, $6111$, and $333$ occurring
respectively $9$, $29$, $1$, and $1$ times.   The remaining $446$ specialization
points yield only $129$ different fields, as for example $(-13/12, 2/9)$, $(11/12, 1/9)$, $(-5/12, -8/27)$,  $(1/4, -1/27)$,  $(1/4, 2/27)$, 
 $(35/108, 8/243)$, $(1/4, 64/3375)$, and $(19/2028, 1/59319)$ all yield the field
 defined by $x^9-9 x^7+27 x^5-27 x^3-4$.  Moreover, a wide variety of subgroups
 of $9T26$ appear, as follows.  
 {\renewcommand{\arraycolsep}{3.8pt}
 \[
 \begin{array}{r|cccccccc}
 \mbox{Group $G$:} & 9T4 & 9T8 & 9T12 & 9T13 & 9T16 & 9T18 & 9T19 & 9T26 \\ 
 \mbox{Size $|G|$:}  & 18  & 36 & 54 & 54 & 72 & 108 & 144 & 432  \\
 \hline
 \mbox{Number of fields in family:} & 2& 1 & 10 & 1 & 5 & 20 & 8 & 82 \\
 \hline
 \mbox{Total number of fields:} & 4 & 1 & 12 & 3 & 5 & 23 & 8 & 87
 \end{array}
 \]
 }%
The last line compares with the relevant complete lists from \cite{JR14}. 
 It gives the total number of number fields with the given Galois group and with discriminant
of the form $-2^a 3^b$ with $a$ even and $b$ odd.    
One can get even a larger fraction of the total number of fields by 
specializing outside of $\scU_{3,1,1}(\Z[1/6])$, both by considering the curve 
 at infinity and then by specializing also at the rare-but-existent
points of say $\scU_{3,1,1}(\Z[1/6p])$, where the auxiliary prime $p$ does 
not divide the discriminant of the field constructed.   The fact that
such a large fraction of all fields of the type considered come from
a single Hurwitz family is suggestive that other Hurwitz families 
may be essentially the only source of number fields with 
certain invariants.  

The current family presents many examples of
phenomena that Principles A, B, and C say are rare in general.  
The drop  from $507$ specialization points to only $129$ fields 
 constitutes many exceptions to Principle A.  
The further drop from $129$ fields 
to just $82$ fields with the generic Galois group 
constitutes many exceptions to Principle B.
  Some of the specializations are tamely ramified or even
unramified at $2$,  and thus correspond to exceptions to Principle C. 
These exceptions form the first data-point arguing for the
expectation already formulated in the introduction: 
 as the complexity of the Hurwitz family increases,
the frequency of exceptions decreases.

\section{A degree 52 family: tame ramification and exceptions to Principle~B}  
\label{deg52}
 In our introductory family, the only exceptions
to Principle B were algebras of the form $K_{h,u} = \Q \times K_{h,u}'$
with $K_{h,u}'$ full.    In specializing many other full families,
most of the exceptions to Principle B we have found have
this very same form. In this section, we present a family which 
is remarkable because some of its 
specializations have a much more 
pronounced drop in fullness. However we do not regard 
this more serious failure of Principle~B 
as anywhere near extreme enough to raise doubts
about Conjecture~\ref{mc}.  

\subsection{A Hurwitz parameter yielding a rational $\AX_h$}
\label{deg52comp}
We start from the normalized Hurwitz parameter 
 \[
h= (S_6,(21111, 222,3_1111, 3_\infty 2_0 1),(2,1,1,1)).
  \] 
  All rational functions with this normalized Hurwitz parameter have
  the form
  \[
  g(s) = \frac{\left(s^3 + b s^2+c s+x\right)^2}{a s^2 (s-y)}.
  \]
  The ramification requirement on $g$ at $1$ is that $(g(1),g'(1),g''(1))=(1,0,0)$.  
  These three equations let one eliminate $a$, $b$, and $c$ via
  \begin{align*}
  a & = -64 (x+1)^2 (y-1)^3, \\   b &= 4 x y-3 x+4 y-6, \\
   c & = -8 x y^2+12 x y-6 x-8 y^2+12 y-3.
  \end{align*}
  Using a resolvent as usual, we find that the critical values of $g(s)$ besides $0$,
  $1$, and $\infty$ are the roots of $W t^2 + (V-U-W) t + U$ where
  \begin{align*}
  U & = (4 x y-x+3 y) \left(64 x^2 y^4-160 x^2 y^3+180 x^2 y^2-108 x^2 y+27 x^2+256 x y^4  \right. \\ & \left. \;\;\;\; -736 x y^3  +864 x y^2-540 x y+162 x+192 y^4-576 y^3+576 y^2-216 y+27\right)^2, \\
  V & = 3^3 (2 x y-x+1)^4 \left(64 x y^3-144 x y^2+108 x y-27 x+64 y^3-144 y^2+81 y\right),\\
  W &= 2^{12} 3^3 (x+1)^4 (y-1)^6 y^3. 
  \end{align*}
 Comparing with the standard quadratic $t^2 + (v-u-1) t + u$, 
one gets the rational presentation
\begin{align}
\label{ratpres52}
u & = \frac{U}{W}, & v & = \frac{V}{W}.
\end{align}
 Summarizing, birationally we have $\AX_{h} = \C_x \times \C_y$,
 $\AU_\nu = \C_u \times \C_v$, and 
 the equations \eqref{ratpres52} give the map $\AX_h \rightarrow \AU_\nu$.  
 Removing $y$ by a resolvent gives the single equation $f_{52}(u,v,x)=0$.   
 Likewise removing $x$ by a resolvent gives the single equation $\phi_{52}(u,v,y)=0$.  
 The left sides have $2781$ and $829$ terms respectively.  
 
 The discriminants of $f_{52}(u,v,x)$ and $\phi_{52}(u,v,y)$ are both
 $-3$ times a square in $\Q(u,v)$.    The Galois groups of these
 polynomials over $\Q(u,v)$ are $S_{52}$, but over $\C(u,v)$
 they reduce to $A_{52}$.  This general phenomenon appeared already in the
 previous section.  It is not of central
 importance to us, which is why we generally refer to 
 full fields and only sometimes make the distinction
 between $S_m$ and $A_m$ fields.

 \subsection{Specialization to curves}  
 \label{specsixcurves}
\begin{figure}[htb]
\includegraphics[width=4in]{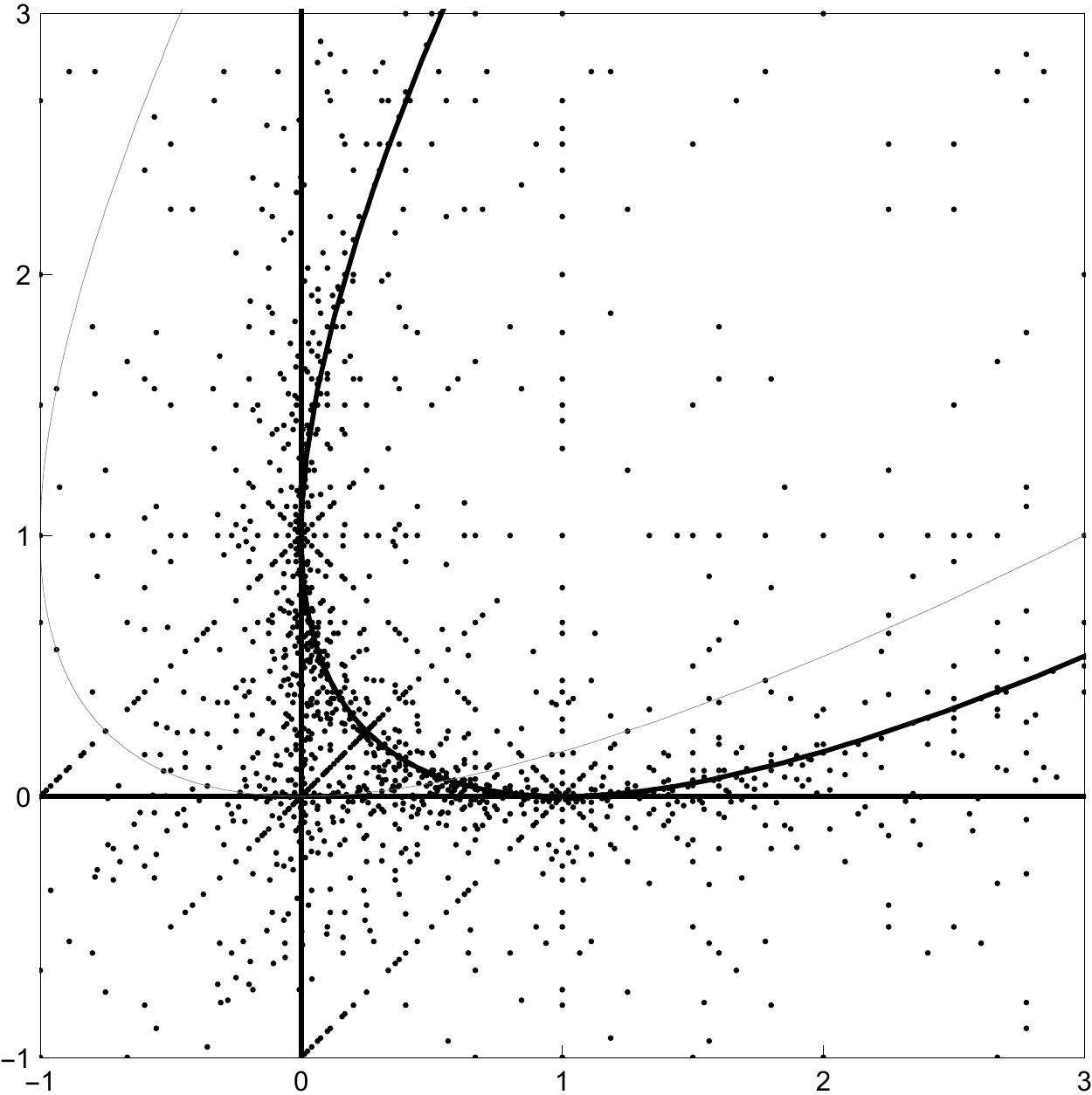}
\caption{\label{spec2111}  A window on $\scU_{2,1,1,1}(\R)$, which is the complement of the three thick discriminantal
curves in the $u$-$v$ plane.  Points are part of the specialization set $\scU_{2,1,1,1}(\Z[1/30])$ which
is used in both \S\ref{deg52} and \S\ref{deg202}.  The thin parabola $(u-v)^2=4 v$ plays a role only in \S\ref{deg52}.  }
 \end{figure}

 Using homogeneous coordinates $U$, $V$, and $W$, related to our 
 standard coordinates $u$ and $v$ via \eqref{ratpres52}, we can
 view $\AU_{2,1,1,1}$ as completed by the projective plane.  Its
 complement in this projective plane has four components, 
 \begin{description}
  \item[$\AAnew$]  the vertical line $U=0$,
  \item[$\AB$]  the horizontal line $V=0$,
  \item[$\AC$] the line at infinity $W=0$, and
  \item[$\AD$] the conic $U^2+V^2+W^2-2 U V - 2 UW - 2 VW=0$.  
  \end{description}
  Figure~\ref{spec2111} draws $\AAnew$, $\AB$, and $\AD$.    
    Note that lines $\AAnew$, $\AB$, and $\AC$ pass through points ${\sf a} = (0:1:1)$, ${\sf b} = (1:0:1)$, and ${\sf c} = (1:1:0)$ respectively, 
  while the conic $D$ goes around ${\sf d} = (1:1:1)$.   
    
 A general line in the projective plane intersects the discriminant locus in five
 points.  However the lines that go through two of the points in $\{{\sf a}, {\sf b}, {\sf c}, {\sf d}\}$ 
 intersect the discriminant locus only three times.  
 These six lines are parametrized in Table~\ref{sixlines}, so that the
 three points become $0$, $1$, and $\infty$.   Exactly as in
 Figure~\ref{scatterintro} earlier, the lines are clearly suggested by the drawn specialization points.   
  Having used homogeneous coordinates for two paragraphs
  to make an $S_3$ symmetry clear, we now return to our standard practice of focusing
  on the affine $u$-$v$ plane.

 When restricted to any one of the six lines, the cover $\AX_h$ remains full.  
 This preserved fullness is in the spirit of Principle B.    
 Table~\ref{sixlines} gives the ramification partitions of these restricted
 covers.  Note that all partitions are even, reflecting
 the fact that the monodromy group is only $A_{52}$.  
 Before beginning any computations with polynomials, 
 we knew these partitions and the fullness
 of the six covers from a braid group computation.  
\begin{table}[htb]
 \[
 \begin{array}{c|cc|lllc}
\mbox{Line} &  u & v & \lambda_0 & \lambda_1 & \lambda_\infty & \mbox{genus}  \\
 \hline
 {\sf ad }&      4t & 1&   12^2 \; 6 \; 4^4 \;  3 \; 2 \; 1 &  3^8 \, 2^{12} \, 1^4 & 20\;  12\;  5 \;   4^3 \,  2 \; 1 & 6  \\
  {\sf bd}  &1 & 4 t &   10^2 \,  8 \;  6^3 \,  5 \; 1  & 3^8 \, 2^{12} \, 1^4 & 10^2 \, 6^2 \,  5 \;  4^2 \,  2^2 \, 1^3 & 5   \\
 {\sf  cd} & t/4 & t/4 &  4^9 \,  2^7 \, 1^2  &  3^8 \, 2^{12} \, 1^4 &  6^4 \, 4^6  \, 2^2 & 0   \\
 {\sf bc }& t & t-1 &     2^{22} \; 1^8 &10^2 \,  8 \;  6^3 \,  5 \; 1 & 6^4 \, 4^6 \, 2^2 & 2   \\
 {\sf ac} & t-1 & t  &    4^9 \, 2^3 \, 1^{10}  &  12^2 \, 6 \; 4^4 \, 3 \; 2 \; 1 & 6^4 \, 4^6 \, 2^2 & 5    \\
 {\sf ab} & t & 1-t &   12^2 \, 6 \; 4^4 \,  3 \; 2 \; 1 & 10^2 \,  8 \;  6^3 \,  5 \; 1 & 5^5 \, 4^2 \, 3^4 \, 1^7 & 9 \\
\end{array}
\]
\caption{\label{sixlines} Six lines in $\AU_{2,1,1,1}$ and topological information on their preimages in
$\AX_{h}$.}
\end{table}

To give an explicit degree $52$ polynomial coming from the cover $\AX_h \rightarrow \AU_{2,1,1,1}$,
we work over the line ${\sf c}{\sf d}$.  The preimage of ${\sf cd}$ is  a curve in the $x$-$y$ plane 
with equation having $x$-degree 3, $y$-degree 6, and twenty-two terms.  A parametrization
is 
\begin{align*}
x & = -\frac{\left(s^2-2 s-2\right) \left(s^4-4 s^3+4 s+2\right)}{2
    s^3}, &
y & = -\frac{3 \left(s^2-2 s-2\right)}{2 (s-2) \left(s^2-4
    s-2\right)}.
\end{align*}
Using the domain coordinate $s$ and the target coordinate $t$, the restricted
rational function takes the form
\begin{equation}
\label{deg52slice}
t = \frac{-A}{C} = \frac{B}{C}+1.
\end{equation}
Here $A$, $B$, and $C$ sum to zero and are given explicitly by
\begin{eqnarray*}
A & = & (s+1)^4\left(s^8-10 s^7+34 s^6-40 s^5-2 s^4+8 s^3+8 s^2+16 s+8\right)^4 \\
    & & \qquad (s-2)^2  \left(s^2-4 s-2\right)^2 \left(s^4-6 s^3+9 s^2-6\right)^2 (s-4) s,  \\
B & = & - \left(s^8-12 s^7+52 s^6-92 s^5+30    s^4+96 s^3-72 s^2-48 s+8\right)^3 \\
     && \qquad \left(s^{12}-12 s^{11}+48 s^{10}-52 s^9-87 s^8+108 s^7+264 s^6-216 s^5 \right. \\
      &&\qquad  \left. -312 s^4+48 s^3+192 s^2+96 s+16\right)^2  \left(s^4-4 s^3+4 s+2\right), \\
C & = & 2^2  \left(2 s^3-9 s^2+6 s+2\right)^6 \left(s^6-6 s^5+6 s^4+10
    s^3-6 s^2-12 s-4\right)^4 \\
    && \qquad  \left(s^2-2 s-2\right)^2.
\end{eqnarray*}
This explicit slice is intended to give a sense of the full cover for $h_{52}$, just 
as the slices in \S\ref{realpictures} and \S\ref{dessin} indicate the
covers for $h_{25}$ and $h^*_{96}$ respectively.   Here we have
explicitly presented information at all three of the cusps, not just
at $0$ and $\infty$ as in  \S\ref{realpictures} and \S\ref{dessin}.   
 
 \subsection{No exceptions to Principles A and C}
Unlike all our previous examples, the current $\nu$ contains at least three ones.  
It thus fits into the framework of \cite{Rob15}, where many $\scU_\nu(\Z[1/\cP])$ 
for such $\nu$ are completely identified.  We therefore can be 
more definitive in reporting specialization results.   
 
The set $\scU_{2,1,1,1}(\Z[1/30])$ contains exactly $2947$ points \cite[\S8.5]{Rob15}
and is drawn in Figure~\ref{spec2111}.  
The Hurwitz number algebras $K_{h,u}$ are all non-isomorphic,
so that Principle A holds without exception.   
All $2947$ algebras are wildly ramified at all three of $2$, $3$, and $5$, 
so that Principle C also holds without exception; in fact 
$\ord_p(D) \geq 52$ fails at $p=2$, $3$, and $5$ only 
$0$, $60$, and $481$ times, so the verification
of Principle C  is particularly easy at $p=2$.

 \subsection{Easily explained exceptions to Principle B} 
 \label{easilyB}
  Twenty-five of the 2947 specialization points
 $(u,v)$ give exceptions to Principle B.  Three of these, namely 
  $(3/8, 1/8)$,  $(1/16, -375/16)$, $(16, -375)$ are exceptions of the 
  sort we have seen earlier: $K_{h,u} = \Q \times K_{h,u}'$ with $K_{h,u}'$ full. 
  Exceptions of this nature are not surprising whenever the Hurwitz cover
  is rational.  In this case the three points in question come respectively
  from points $(-1/6,3/8)$, $(-4/3,3/4)$, and $(-3,3/4)$ in $\scX_h(\Q)$. 
  
  In fact, any point  $(x,y) \in \scX_h(\Q)$ causes such a factorization, but typically
  its image $(u,v)$ causes ramification at extraneous primes.   
  For example, take $s=1$ in \eqref{deg52slice}, making 
  \begin{equation}
  \label{deg51abc}
  t = \frac{111936400}{43923} = \frac{2^4 \,  5^2 \,  23^4}{3 \; 11^4} = 1 + \frac{37^3 \, 47^2}{3 \; 11^4}.
  \end{equation}
  The number field $K'$ defined by the degree $51$ factor of $f_{52}(t,s)$ has discriminant 
  $-2^{102} 3^{95} 5^{48} 11^{16} 37^{12} 47^{16}$.  Whenever we discuss exceptions to Principle B,
  we always have in mind a fixed $\cP$, here $\{2,3,5\}$, and do not consider fields
  like $K'$ to be exceptions.  
  
  \subsection{Ramification at tame primes} 
  \label{ramificationtame}
  Recall from \S\ref{wildram} that ramification at $p$ in a Hurwitz number algebra 
  $K_{h,u}$ can only be wild if $p \in \cP_h$ or $p \leq \max_i \nu_i$.   The field $K'$ from
  the previous subsection presents a convenient opportunity to illustrate how 
  ramification in $K_{h,u}$ at the remaining primes is calculable in purely 
  group-theoretic terms.  
  
  To describe the factorization of the local algebras $K_p'$, we represent the fields 
  appearing by symbols $e_c^f$ as in \eqref{fieldsymbols}. 
   We simplify by just writing $e^f$ for tame 
  fields, since tameness implies $c=e-1$.  The factorizations are
   \begin{align*}
  \mbox{2:} & \; {\bf 16}_{38}  \, {\bf 16}_{38} \, {\bf 8}_{16} \,  {\bf 2}_3 \, {\bf 2}_3 \, {\bf 2}_2 \,  {\bf 2}_2  \, 1^2  \,1, \;\;\;\;\;\;\;\;\;\;&
  \mbox{11:} & \; 3^6 \, 3^2 \, 1^{16} \,  1^2 \, 1^2 \,  1^2 \, 1^2 \,  1 \; 1 \; 1 &\!\!\!\!\!& \rightarrow 3^8 1^{27}, \\
  \mbox{3:}   & \; {\bf 18}_{39} \, {\bf 12}_{21}^2 \, {\bf 6}_{11} \; {\bf 3}_3, &
  \mbox{37:}  & \; 2^5 \, 2^4 \,  2^2 \,  2 \;1^{18} \, 1^3 \, 1 \; 1 \; 1 \; 1 \; 1 \; 1 &\!\!\!\!\!& \rightarrow 2^{12} 1^{27}, \\
  \mbox{5:}  & \; {\bf 25}_{40} \, 2^4 \, 2^4 \, 1^4 \, 1^2  \, 1 \; 1 \; 1 \; 1, &
  \mbox{47:}  & \; 3^6 \, 3^2 \,  1^4 \,  1^4 \, 1^4 \, 1^4  \,1^4 \, 1^2 \, 1 &\!\!\!\!\!& \rightarrow 3^8 1^{27}.
  \end{align*}
  The wild primes behave in a complicated way as always, with $p$-wildness at $p=2$, $3$, and $5$
  being $48$, $51$, and $25$.  However the tame primes are much 
  more simply behaved.  
  
  To work at an even simpler level, we factor over the maximal unramified extension
  of $\Q_p$, rather than $\Q_p$ itself.  For tame primes, this corresponds to regarding
  the printed exponents $f$ simply as multiplicities, and collecting together symbols with
  a common base.   The resulting tame ramification 
  partitions are indicated to the right, after arrows.  
  
  Note that there are actually four primes greater than $5$ involved in \eqref{deg51abc}.  With
  their naturally occurring exponents, $11^4$ is associated to $\infty$, 
  $23^4$ to $0$, and $37^3$ and $47^2$ to $1$.   In general, tame ramification partitions
  can be computed from the placement of the specialization point in $\scU_\nu(\Q_p)$
   and braid group considerations.  In the setting of three-point covers, the general formula is simple, and uses
   the standard notion of the power of a partition.  Namely,  if $p^m$ is associated to $\tau \in \{0,1,\infty\}$
   its tame ramification is the power $\lambda^m_\tau$ of the
   geometric ramification partition $\lambda_\tau$.   Applying the ${\sf cd}$ line of Table~\ref{sixlines},
   the partitions $\lambda_\infty^4=3^8 1^{27}$, $\lambda_0^4 = 1^{51}$, $\lambda_1^3 = 2^{12} 1^{27}$,
   and $\lambda_1^2 = 3^8 1^{27}$ do indeed agree with the partitions found by 
   direct factorization of the polynomial defining $K'$.     
   
   The mass heuristic reviewed in \S\ref{ffmh} is based on an equidistribution principle.  In the horizontal direction, it translates to the following conjecture, proved for $m \leq 5$:
    when one considers full degree $m$ fields ordered by their absolute discriminant outside of $p$, all tame ramification partitions are asymptotically equally likely.  
    We regard the fact that Hurwitz number fields
   escape the mass heuristic as being directly related to their highly structured ramification.  
   In the current instance, there are $239,943$ partitions of the integer $51$, and the two partitions 
   $2^{12} 1^{27}$ and $3^8 1^{27}$ are far from typical.        
  
   \subsection{A curve of more extreme exceptions to Principle B}
 The twenty-two exceptions not discussed in \S\ref{easilyB} all have a common geometric source: 
 above the base curve $\SB$ given by $(u-v)^2=4v$, the cover splits
 into a full degree $42$ cover of genus five and a full degree $10$ 
 cover $\SC$ of genus zero.     While decompositions $52=51+1$ are
 governed by rational points on $\scX_h$ itself, decompositions 
 of the form $52=42+10$ are governed by rational points on a
 resolvent variety of degree $\binom{52}{10}$ over
 $\scU_{2,1,1,1}$.  As this degree is about 15 billion, 
 the existence of a entire curve of rational points is remarkable.

 To reveal the structure of the cover $\SC \rightarrow \SB$,
 we parametrize the base curve $\SB$ via
 \begin{align}
 \label{parabparam}
 u & = \frac{4t}{(t-1)^2}, & v & = \frac{4}{(t-1)^2}.
 \end{align}
   In the decompositions $K_{h,(u,v)} = K^{42}_{t} \times K_t$,
 the twenty-two  $K^{42}_t$ are all full degree forty-two fields, with pairwise distinct discriminants.   
 
    The genus zero curve $\SC$ is given by $x (4 y-3)^3=-24 y^2 (2 y-3)$, and 
    so $y$ is a parameter.   The map from the $y$-line $\SC$ to 
    the $t$-line $\SB$ is given by the vanishing of
    \begin{eqnarray*}
    f_{10}(t,y) & = & (4 y-3) (8 y-3) \left(32 y^4-192 y^3+360 y^2-252 y+27\right)^2 \\ && \qquad  \; \,+ t (4 y-9) \left(96 y^4-256
    y^3+216 y^2-108 y+27\right)^2.
    \end{eqnarray*}
    Thus one has two visible ramification partitions $\lambda_0 = \lambda_\infty = 222211$.  The discriminant
    of $f_{10}(t,y)$ is $-2^{136} 3^{57} 5^{25} t^4 (t-1)^5 (t-9)^5$.  At the other singular values,
    the ramification partitions are $\lambda_1 = \lambda_9 = 32221$.    In 
    fact, the decic algebras $K_t$ and $K_{9/t}$ are isomorphic via the involution $y \mapsto (6y-9)/(8y-6)$.

At the level of the decic cover only, we have just indicated a failure of Principle A: rather
than $22$ distinct decic algebras, there are ten pairs switched by $t \leftrightarrow 9/t$ and
then two algebras $K_3$ and $K_{-3}$ arising once each.   The ten algebras arising twice are all
full fields and wildly ramified at all three of $2$, $3$, and $5$.   However $K_{3}$ and
$K_{-3}$ are not full, and not wildly ramified at $5$, giving 
failures of Principle B and C at this decic cover level.  

In terms of supporting Conjecture~\ref{mc} for $\cP = \{2,3,5\}$, the exceptional 
behavior above $\SB$ is in a sense good.   Instead of twenty-two 
contributions to $F_{\cP}(52)$, one gets twenty-two contributions
to $F_{\cP}(42)$ and then ten more to $F_{\cP}(10)$.   But 
in another sense this exceptional behavior is bad.  It explicitly
illustrates phenomena which, if occurring ubiquitously
in high degree, might make Conjecture~\ref{mc} false. 
However our computations suggest that, far from becoming
ubiquitous, the phenomena exhibited here become rarer as degrees
increase.

\section{A degree $60$ family: non-full monodromy and a prime drop}
\label{deg60}
     The statement of Conjecture~\ref{mc} involves all finite nonabelian
 simple groups equally.   In this paper, however, we are focusing  
 on the simple groups $A_5$ and $A_6$ because of the computational
 accessibility of the corresponding families $\AX_h \rightarrow \AU_\nu$.
 In this section and the next, we add some balance
 by  presenting results on covers coming from simple groups 
 not of the form $A_n$.   The family presented here has
 the particular interest that it is non-generic in two ways.  
  
 \subsection{A Hurwitz parameter with unexpectedly non-full monodromy} 
 \label{konigbasic}
 The simple group $G = PSL_3(\F_3)$ has order $5616 = 2^4 \cdot  3^3 \cdot 13$ and outer automorphism
group of order two.    It has two non-isomorphic degree $13$ transitive permutation representations, coming from an action on a projective plane $\bbP^2(\F_3)$ and
its dual $\widehat{\bbP}^2(\F_3)$.   These actions are interchanged by the outer
involution.  The two smallest non-identity conjugacy classes in $G$ consist
of order $2$ and order $3$ elements.  
In each of the degree $13$ permutation representations,
these elements act with cycle structure $2^4 1^5$ and $3^3 1^4$ respectively.  

Let $h = (PSL_3(\F_3),(2^4 1^5, 3^3 1^4),(3,2))$.  To conform to our main reference \cite{Kon} for this
section, we make a quadratic base change and work over
 $\AU_{3,1,1}$ rather than $\AU_{3,2}$.  A braid group computation 
reveals that the degree $120$ 
Hurwitz cover $\pi_h$ factors as a composition of three covers as indicated:
\begin{equation}
\label{threecovers}
\AX_h \stackrel{2}{\rightarrow} \AX^*_h \stackrel{15}{\rightarrow} \AQuart \stackrel{4}{\rightarrow}  \AU_{3,1,1}.
\end{equation}
The intermediate cover $\AX^*_h$ is just the quotient of $\AX_h$ by the natural action of $\Out(G)$.   
This failure of fullness illustrates one of the general phenomena treated at length in \cite{RV15}.  

However, very unusually in comparison with Table~\ref{accessible}, 
the reduced Hurwitz cover
$\pi_h^* : \AX^*_h \rightarrow \AU_{3,1,1}$ is also not full.  It clearly fails
to be primitive, because of the intermediate cover $\AQuart$.  Moreover, the  degree fifteen
map is not even full, as its monodromy group is $S_6$ in a degree $15$  transitive 
representation. 

The degree $13$ covers of the projective line parametrized by $\AX_h$ have
genus zero.   Using this fact as a starting point, K\"onig \cite{Kon} succeeded in finding 
coordinates $a$, $b$ on $\AX_h$, with corresponding covers $\AP^1_s \rightarrow \AP^1_t$ being as follows.    Define 
\begin{eqnarray*}
f_0 & = & \frac{a b s}{3}+\frac{a b}{9}+a s^2-\frac{a}{3}+s^3, \\
f_1 & = & \frac{s^2 \left(a b^2-4 a b+12 a-3 b^2-9\right)}{(b-3)^2}+\frac{s \left(a
   b^2-4 a b+12 a-9 b-9\right)}{3 (b-3)}+s^3-1, \\
g_0 & = & \frac{a b s}{3}+\frac{a b}{9}+a s^2-\frac{a}{3}+s^3, \\
g_1 & = & \frac{1}{9} s \left(4 a b^2-6 a b+9 a+9 b-27\right)+\frac{1}{3} s^2 (4 a
   b-3 a+9)+a s^3-a.
\end{eqnarray*}
Then the two-parameter family is given by
$g(a,b,t,s) := f_0^3 f_1 s - t g_0^3 g_1 =  0$.
 
K\"onig's interest in this family is in producing number fields 
with Galois group $G$.  For example $(a,b,t) = (-9,-6,-3)$ gives
a totally real such field with discriminant $3^{12} 251^4 353^4$.  
To systematically study specializations, it is important to determine
the discriminant of $g(a,b,t,s)$.   Computation shows
that it has the
following form:
\begin{align*}
D(a,b,t) & = 
                \left(-\frac{4}{3} a b^3 +a^2 b^2+6 a b^2-3 b^2-4 a^2 b-18 a
    b+18 b+12 a^2-27\right)^{28} \\
       &  \qquad  \qquad a^{12} (b-3)^{18} t^6 \left(C_0 t^3 + C_1 t^2 + C_2 t + C_3 \right)^4.
\end{align*}
Here $C_0$, $C_1$, $C_2$, and $C_3$ as expanded elements
of $\Q[a,b]$ have $24$, $45$, $53$, and $36$ terms respectively. 
Because of the complicated nature of this discriminant, it 
is hard to get field discriminants to be as small as the one
exhibited above.  

For K\"onig's purposes of constructing degree thirteen fields with Galois group $G$, he does not 
need the map to configuration space at all.  To move over into
our context of constructing Hurwitz number fields, we do need this map.
Replacing $t$ in $(C_0 t^3 + C_1 t^2 + C_2 t + C_3)$ with $C_1 t/C_0$ and
setting the resulting cubic proportional to $t^3 + t^2 + u t + v$ gives a degree $120$ map 
$\pi_h$ from the $a$-$b$ plane $\AX_h$ to the $u$-$v$ plane $\AU_{3,1,1}$.   Removing
$a$ from the pair of equations gives a degree $120$ polynomial $f_{120}(u,v,b) \in \Q(u,v)[b]$ 
describing the covering map.       

\subsection{Reduction to degree $60$}    To reduce from the degree $120$ cover $\AX_h$ to the
degree $60$ cover $\AX_h^*$, we proceed as follows.  For $(a_i,b_i) \in \Q^2$,
one gets $(u_i,v_i) = \pi_h(a_i,b_i) \in \Q^2$.  Then $f_{120}(u_i,v_i,b) \in \Q[b]$ factors.
For almost all choices of $(a_i,b_i)$, the degrees of the irreducible factors 
are $90$, $6$, $6$, $4$, $4$, $4$, $4$, $1$, and $1$.   One of the
linear factors is $b-b_i$ and we write the other one as $b-b_i'$.  Then
typically just one rational number $a_i'$ satisfies the two equations
$\pi_h(a'_i,b_i') = (u_i,v_i)$.    From enough datapoints we interpolate to 
get the canonical involution on $\AX_h$.  It is
\begin{align}
\label{koniginvol}
a' & = \frac{(b-3) (4 a b+6 a+9)}{a b^2-4 a b+12 a-18 b+18}, &
b' & = \frac{3 b}{b-3}.
\end{align}
This involution is useful even in K\"onig's context.  For example,
specializing at $(a',b',t) = (171/58,2,-3)$ gives the 
dual totally real number field, also with discriminant $3^{12} 251^4 353^4$.

A quantity stabilized by the involution is $x = b^2/(b-3)$.   The resolvent 
$$\mbox{Res}_b(f_{120}(u,v,b),(b-3) x-b^2)$$
is proportional the square of degree $60$ polynomial $f_{60}(u,v,x)$.   This polynomial captures
the cover $\AX_h^* \rightarrow \AU_{3,1,1}$.  

\subsection{Low degree resolvents} From the braid group 
computation, we know that the monodromy group has quotients
of type $S_3$, $S_4$, and $8T40 = 2^3.S_4$.  Here the 
$S_4$ quotient corresponds to the cover $\AQuart$.  Equations
for these quotients and their discriminants are  
\begin{align*}
f_3(u,v,x) & = x^3+x^2+x u+v,  & D_3 & =  d, \\
f_4(u,v,x) & = x^4-2 x^2 v-8 x v^2-4 u v^2+v^2, & D_4 & = - 2^{12} d v^6, \\
f_8(u,v,x) & = x^8+8 x^4 d u v-72 x^4 d v^2 & D_8 & = -2^{60} d^{17}  
v^{14} \left(u^3-v\right)^4. \\
&  \qquad +64 x^2 d^2 v^2-16 d^3 v^2,
\end{align*}
Here we have seen the cubic and quartic polynomials in \S\ref{normalization}, with
$d$ being given explicitly in \eqref{quartdisc}.

\subsection{Reduction to degree $24$}  The equation $f_4(u,v,m)=0$ 
is linear in $u$.   Solving it gives  $u = (m^4 - 2 m^2 v + v^2 - 8 m v^2)/(4 v^2)$.
Expressing $f_{60}(u,v,x)$ in terms of $m$, $v$, and $x$ and factoring, one gets
$g_{15}(m,v,x) g_{45}(m,v,x)$.  Here $g_{15}(m,v,x)$ has Galois group $S_{6}$
over $\Q(m,v)$, in a degree $15$ permutation representation.  

Abbreviate $e=m^3-m v-2 v^2$. Then the polynomial for the
standard sextic representation works out to 
\[
g_6(m,v,x)  =  2 x^6 v^2-3 x^4 e \left(m^2-v\right)-8 x^3 e^2-6 x^2 e^2 m+2 e^3.
\]
Returning to the original base, one gets a degree $24$ polynomial,
\[
f_{24}(u,v,x) = \mbox{Res}_m(f_4(u,v,m),g_6(m,v,x)).
\]
Similarly, by means of the outer automorphism of $S_6$, one 
has a twin polynomial $g_6^t(m,v,x)$ and its degree $24$ polynomial 
$f^t_{24}(u,v,x)$.  While $f_{60}(u,v,x)$, $f_{24}(u,v,x)$, and $f_{24}^t(u,v,x)$
all have the same splitting field, the latter two are much easier to 
work with because of their lower degree.  

\subsection{Specialization to number fields} 
We have specialized at the $507$ points in $\scU_{3,1,1}(\Z[1/6])$ considered in \S\ref{deg9},
obtaining $507$ algebras with discriminant of the form $\pm 2^a 3^b$.   
Replacing $(u,v)$ by $(v/u^2,v^2/u^3)$, corresponding to the involution 
of $\AU_{3,1,1}$
with quotient $\AU_{3,2}$, gives an isomorphic algebra.  
We report on the fields involved in these algebras,
 since Galois groups are small enough so that future comparison
with other sources of fields with these groups seems promising.  

For simplicity, we exclude the twenty-three $(u,v)$ where $u=0$, so that the involution 
above is everywhere defined.  We switch coordinates to the coordinates used in \cite{RobG2} via 
$(p,q) = (3u,3v/u^2)$ and $(u,v) = (p/3,p^2q/27)$.  
In the new coordinates, the involution is simply $(p,q) \mapsto (q,p)$, and we normalize
by requiring $p \leq q$.   We then have $232$ algebras $K_{p,q}$ with $p < q$ and 
$20$ algebras $K_{p,p}$.  Besides these algebras, we have their twins $K_{p,q}^t$,
and their common octic and quartic resolvents $\tilde{R}_{p,q}$ and $R_{p,q}$.

Despite the non-generic behavior of the family in general, Principal A has
no exceptions in the current context: the $252$ 
algebras $K_{p,q}$ and their $252$ twins $K_{p,q}^t$ form $504$ 
non-isomorphic algebras.  Principle C also has no exceptions, as all 
algebras are wildly ramified at both $2$ and $3$.  

There are many exceptions to Principle B.  For example $K_{153/1849,129/289}$ 
factors as $6+6+12$ with the factors having Galois group $6T9$, $6T15=A_6$, 
and $12T299 = S_6 \wr S_2$.  Its twin factors as $3+3+6+12$ with factors
having Galois groups $S_3$, $S_3$, $A_6$, 
and $S_6 \wr S_2$.   The two $A_6$ factors are given by the polynomials
\begin{eqnarray}
\label{smalla6} f_6(x) & = & x^6 - 3 x^5 + 3 x^4 - 6 x^2 + 6x - 2, \\
\label{smalla6t} f^t_6(x) & = & x^6 - 3 x^4 - 12 x^3 - 9 x^2 + 1.
\end{eqnarray} 
These polynomials will be discussed further at the end of the next subsection.  

For the rest of this section, we avoid Galois-theoretic complications
like those of the last paragraph by 
requiring that $R_{p,q}$ either has an irreducible cubic 
factor or is irreducible itself.   There are 39 $(p,q)$ of the
first type, and 178 $(p,q)$ of the second.  
Failures of Principle B in this restricted setting
are very mild, as $A_6^4$ is a subgroup of the Galois group of all these specializations.  
In the case of a cubic-times-linear quartic resolvent, we change 
notation by focusing on the larger degree part, so that $K_{p,q}$, $K_{p,q}^t$, $\tilde{R}_{p,q}$,
and $R_{p,q}$ now have degrees $18$, $18$, $6$, and $3$.   

\subsection{Some lightly ramified number fields}
Table~\ref{konigspecs} summarizes the fields under consideration, with 
resolvent Galois groups indicated by $Q$ and $\tilde{Q}$.    
In all cases, if $K_{p,q}$ has some Galois group $mTj$ then its
twin $K_{p,q}^t$ has the same Galois group $mTj$.  
For each Galois group,
the table gives a corresponding field in our collection with smallest
root discriminant.  Thus $(p,q)$ is chosen because one of $\delta = \mbox{rd}(K_{p,q})$ and
$\delta^t = \mbox{rd}(K_{p,q}^t)$ is small; the other is sometimes substantially 
larger.    Galois groups were computed by {\em Magma}, making use thereby
of the algorithms of \cite{FK14} and works classifying permutation groups.  

\begin{table}[htb]
\[
{\renewcommand{\arraycolsep}{2pt}
\begin{array}{crc|cc|cr|rc|rc}
Q & |\tilde{Q}| & \tilde{Q} & p & q & \mbox{Gal}(K_{p,q}) & \# & \multicolumn{1}{c}{D} & \delta &  \multicolumn{1}{c}{D^t} & \delta^t \\
\hline
S_3 & 6 & 6T2 & 1/12 & 24 & 18T971 &1& 2^{22} 3^{46} & 38.66 & 2^{24} 3^{44} & 36.95 \\ 
S_3 & 12 & 6T3 & -3/125 & 1 & 18T972&7 & 2^{26} 3^{41}& 33.24 & 2^{26} 3^{39} & 29.42 \\ 
A_3 & 6 & 6T6 & 9/121 & 11/27 & 18T974 &3& - 2^{18} 3^{46} & 33.14  & -2^{24} 3^{44} & 36.95 \\
S_3 & 24 & 6T7 & -3/49 & 7/3 & 18T976 &1& 2^{30} 3^{46} & 52.60 & 2^{32} 3^{44} & 50.29 \\ 
S_3 & 48 & 6T11 & 1/9 & 9 & 18T977 &27& -2^{33} 3^{41} & 43.52 & -2^{27} 3^{39} & 30.57 \\
\hline
C_4 & 16 & 8T7 &  {18}/{25} & {5}/{6} & 24T24946 &2
& 2^{87} 3^{24} & 37.01 & 2^{97} 3^{24} & 49.41  \\
V & 16 & 8T9 & -72 & {3}/{16} & 24T24948 &3
& 2^{54} 3^{40} & 29.68 & 2^{70} 3^{38} & 43.00   \\ 
V & 16 & 8T11 &  -27 & -{1}/{3} & 24T24949 &3
& 2^{48} 3^{38} & 22.78 & 2^{52} 3^{42} & 30.70  \\
D_4 & 16& 8T6 &  -9 & {1}/{3} & 24T24952 &8
& 2^{81} 3^{24} & 31.12 & 2^{85} 3^{24} & 34.94 \\
D_4 & 16 & 8T8 & 1 & {9}/{8} & 24T24953 &5
& -2^{52} 3^{43} & 32.14 & - 2^{60} 3^{45}  & 44.38   \\
A_4 & 24 & 8T13 & -2 & {1}/{4} & 24T24956 &5
& 2^{54} 3^{46} & 39.07 & 2^{56} 3^{46}  & 41.39  \\
D_4 & 32 & 8T17 &  -{8}/{3} & {9}/{16} & 24T24961 &3
& 2^{83} 3^{24} & 32.97 & 2^{81} 3^{24} & 31.12   \\
D_4 &32 & 8T15 & -{27}/{25} & {5}/{9} & 24T24962 &17
& 2^{52} 3^{43} & 32.14 & 2^{48} 3^{45} & 31.38      \\
C_4 & 32 & 8T16 &  {16}/{27} & {27}/{32} & 24T24964 &1
& 2^{79} 3^{24} & 29.38 & 2^{97} 3^{24} & 49.41  \\
S_4 & 48 & 8T23 &  -{1}/{3} & 1 & 24T24968 &44
& -2^{46} 3^{59} & 56.22 & - 2^{46} 3^{51} & 38.98  \\
D_4 & 64 & 8T26 & -{135}/{289} & {17}/{25} & 24T24974 &26
& 2^{58} 3^{35} & 26.50 & 2^{50} 3^{39} & 25.26   \\
S_4 & 192 & 8T40 & -{7}/{12} & {32}/{49}  & 24T24982 &61
& 2^{55} 3^{59} & 72.91 & 2^{49} 3^{51} & 42.51\\
\end{array}
}
\]
\caption{\label{konigspecs} Fields $K_{p,q}$ and $K_{p,q}^t$ 
with given Galois group and small root discriminant}
\end{table}

For almost all groups in degree $\leq 19$, the database of Klueners 
and Malle presents at least one corresponding  field.  The database also highlights
the field presented with smallest absolute discriminant.   
For the five degree eighteen groups appearing in Table~\ref{konigspecs},
our fields are well under the previous minima, these
being  $643.84$, $51.78$, $66.63$, $71.35$, and $57.52$ in the order listed.   
 For the twelve degree twenty-four groups, we similarly
do not know of other fields with smaller root discriminants.  

The light ramification in these fields is often reflected in the smallness of coefficients in the
standardized
polynomials returned by {\em Pari}'s \verb@polredabs@. 
For example, the degree eighteen field in the table of smallest root discriminant is $K_{-3/125,1}^t$.  
It is defined by 
\begin{eqnarray*}
\lefteqn{f_{18}(x) =} \\
&& x^{18}+9 x^{16}-18 x^{15}+18 x^{14}-36 x^{13}+72 x^{12}-18 x^{11}+36
    x^{10} \\ && -180 x^9+18 x^8+54 x^7+48 x^6-108 x^5+18 x^4-30 x^3+9 x^2-1. 
\end{eqnarray*}
Similarly, the degree twenty-four field in the table of smallest root discriminant is $K_{-27,-1/3}$.  
It is defined by 
\begin{eqnarray*}
\lefteqn{f_{24}(x) = x^{24}-8 x^{21}+64 x^{18}-36 x^{17} } \\
&& -9 x^{16}-56 x^{15}+276 x^{14}-72
    x^{13}+237 x^{12} -24 x^{11}  +486 x^{10}-88 x^9 \\ && +513 x^8+36 x^7+256 x^6+48
    x^5+18 x^4+20 x^3-6 x^2+1.
 \end{eqnarray*}
 Another particularly interesting case comes from the second to last line of Table~\ref{konigspecs},
 where both $\delta$ and $\delta^t$ are small.   
 
 For speculating where the fields of this section may fit into complete lists, it is
 insightful to compare with the polynomials from \eqref{smalla6} and \eqref{smalla6t}.  
 The fields $\Q[x]/f_6(x)$ and $\Q[x]/f^t_{6}(x)$ have 
 root discriminants $\delta = (2^8 3^8)^{1/6} \approx 10.90$ and $\delta^t = (2^{10} 3^8)^{1/6} \approx 13.74$.
 These root discriminants are $12^{\rm th}$ and $44^{\rm th}$ on the complete sextic $A_6$ list,
 substantially behind the first entry $(2^6 67^2)^{1/6} \approx 8.12$ \cite{JR14}.  On the other
 hand the common splitting field of $f_6(x)$ and $f^t_6(x)$ has root discriminant 
 $2^{13/6} 3^{16/9} \approx 31.66$.  This is the smallest root discriminant of a Galois $A_6$ field,
 substantially ahead of the second smallest $2^{7/6} 3^{25/18} 13^{1/2} \approx 37.23$ \cite{JR14}.  
  We expect that the degree $18$ and $24$ fields discussed in this subsection behave similarly
  to these sextic fields:
  their root discriminants should appear early on complete lists, and their
  Galois root discriminants should appear even earlier.

\section{A degree $96$ family: a large degree dessin and Newton polygons}
\label{deg96}  Almost all the full number fields presented so far
in this paper have been ramified 
exactly at the set $\{2,3,5\}$.   Conjecture~\ref{mc} on the other hand envisions
inexhaustible supplies of full number fields with other ramification sets 
$\cP$.  This section presents examples with $\cP = \{2,3,7\}$.  

\subsection{One of two similar Hurwitz parameters}
Theorems 4.1 and 4.2 of Malle's paper  \cite{Mal00} each give a two-parameter family of septic 
covers of the projective line with monodromy group $SL_3(\F_2) \subset S_7$ of order
$168 = 2^3 \cdot 3 \cdot 7$.  We focus on the family of Theorem~4.2 which is indexed by
the Hurwitz parameter
\[
h = (SL_3(\F_2), (22111,421), (4,1)).
\]
The corresponding degree is $m = 192$.

The situation has much in common with K\"onig's situation from \S\ref{konigbasic}
and we can proceed similarly.  
 Thus, by equating a discriminantal factor with a standard quartic,
we realize $\AX_h$ as a degree $192$ cover of 
$\AU_{4,1}$.   The outer involution of $SL_3(\F_2)$ coming from projective duality 
gives an explicit involution analogous to \eqref{koniginvol}.  Quotienting by this involution
yields the degree $96$ cover 
$\AX^*_h \rightarrow \AU_{4,1}$.   
Unlike the exceptional cover of the previous section, this cover is full.   

The family from Theorem~{4.1} of \cite{Mal00} is very similar: the partition $421$ is replaced by $331$, 
and the degree $192/2=96$ is replaced by $216/2=108$.   We 
are working with the degree 96 family because the curve given by $f_{96}(j,x)=0$  below 
has genus zero, while its analog for the degree $108$ family
has genus one.   

\subsection{A dessin} 
\label{dessin}
The reduced configuration space $\scU_{4,1}$ is the same as that for 
our introductory family and has been described in \S\ref{normalization}.  
However the specialization set is now 
 $\scU_{4,1}(\Z[1/42])$ rather than the
$\scU_{4,1}(\Z[1/30])$  drawn in Figure~\ref{scatterintro}.
We present here only a polynomial for the degree 
$96$ cover of the vertical
line $(u,v)=(1/3,(j-1)/27 j)$ evident in Figure~\ref{scatterintro}:
{\small
\begin{eqnarray*}
\lefteqn{f_{96}(j,x) =} \\
  & & \left(7411887 x^{32}-316240512 x^{31}+5718682592 x^{30}-57608479936
    x^{29} \right. \\ && +345466405984 x^{28}-1143902168192 x^{27}+500924971008
    x^{26}+20121596404224 x^{25} \\ && -178485128485440 x^{24}+1076315934382080
    x^{23}-4902849972088320 x^{22} \\ &&
    +16964516971136000 x^{21}-45252388465854976
    x^{20}+95197078307043328 x^{19} \\ && -161987009378324480
    x^{18}+229049096903122944 x^{17}-277106243726667264
    x^{16} \\ && +295558502345637888 x^{15}-284898502452436992
    x^{14}+250987121290100736 x^{13} \\ && -200876992270295040
    x^{12}+143474999551229952 x^{11}-89556680876359680
    x^{10} \\ && +47950288840949760 x^9-21681369027919872 x^8+8162827596988416
    x^7 \\ && -2520589064601600 x^6+626540088655872 x^5-122178152300544
    x^4 \\ && \left. +17986994307072 x^3-1878160048128 x^2+123834728448 x-3869835264\right)^3 \\
    && -2^{20} j x^6 (3 x-2)^2 \left(x^2+2 x-2\right)^6 \left(7 x^2-14
    x+6\right)^{21} \left(2 x^3-15 x^2+18 x-6\right)^9.
\end{eqnarray*} 
}

\noindent
The printed degree thirty-two polynomial
capturing behavior at $j=0$ has Galois group $A_{32}$ and
 field discriminant only $2^{64} \, 3^{36} \, 7^{18}$.

\begin{figure}[htb]
\includegraphics[width=4.6in]{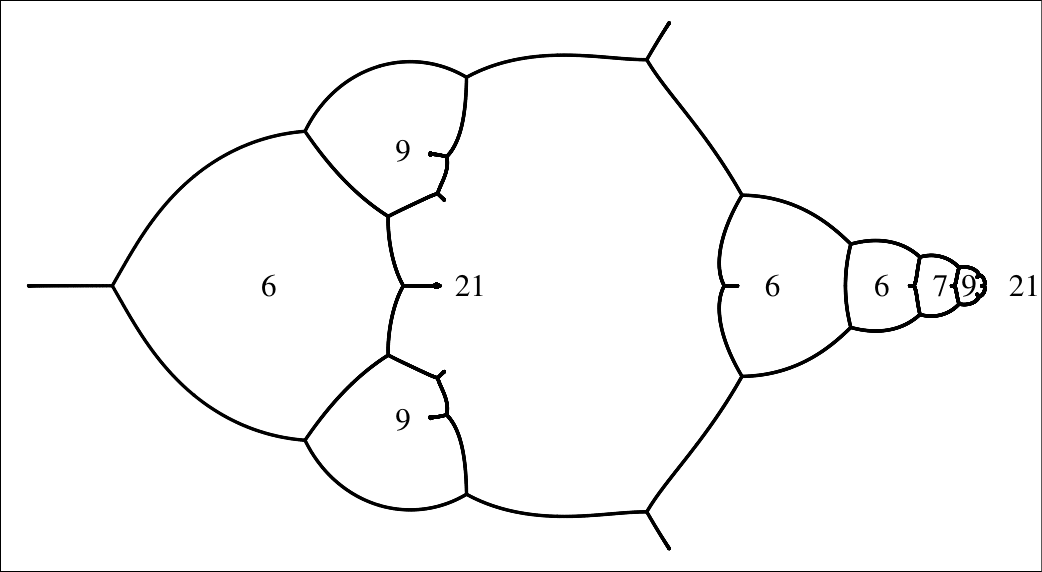}
\caption{\label{dessin96} The dessin corresponding to $f_{96}(j,x)$.  Besides the nine 
regions with indicated sizes there is a tenth region of size $1$ immediately to the left
of the centrally printed $21$.  This small region is adjacent to two triple points and near an endpoint.  
Also to the immediate right of each of the two left-printed nines, there is one triple point
and two endpoints. }
\end{figure}

Figure~\ref{dessin96} draws the dessin of $f_{96}(j,x)$, not in the copy of
$\C$ with coordinate $x$, but rather the copy of $\C$ with coordinate 
$x'=1/(1-x)$, for better geometric appearance.   By definition, the figure consists
of all $x'$ corresponding to $x$ satisfying $f_{96}(j,x)=0$ with $j \in [0,1]$.
This figure has the natural structure of a graph with $96$ edges, the preimages
of $(0,1)$.  All vertices have degree $\leq 3$:  there are thirty-two triple points,
the preimages of $0$, and forty double points and sixteen endpoints,
the preimages of $1$.   The forty double points are not readily visible in
the figure, as they lie in the middle of forty double edges, but most of the triple 
points and endpoints are.   There are also
ten regions, of varying size, defined as half the number of bounding edges.
The few aspects of all this structure which are not visible are described
in the caption of Figure~\ref{dessin96}.    The topological structure
could also be deduced from a braid computation, rather
than from the defining equation.  

    The polynomial $f_{96}(j,x)$ and Figure~\ref{dessin96} 
 illustrate the nature and complexity of the
 objects we are considering.   Note that the existence of this cover
 shows that the Hurwitz number algebra indexed by $(A_{96}, \; (3^{32}, \, 2^{40} \, 1^{16}, \, 21^2 \,  9^3 \, 7 \;  6^3 \, 2), \; (1,1,1))$ has at least one factor of $\Q$.   The entire Hurwitz
 algebra is way out of computational range, because the two main 
 terms in the mass formula \eqref{massformula}
 give $3 \times 10^{15}$ as an approximation for its degree.  
 
    A common feature of $f_{25}(j,x)$ from \S\ref{realpictures} and $f_{96}(j,x)$ 
  is not accidental.  In the braid group description of their monodromy,
 calculable purely group-theoretically, local monodromy operators about $0$ and
 $1$ are the images of braid group elements of order $3$ and $2$ respectively.
 Thus the preimage of $u=1/3$ in $\AX_h \rightarrow \AU_{4,1}$ for
 any $h$ with multiplicity vector $(4,1)$ likewise has this property.

 \subsection{Specialization and Newton polygons} 
  For greater explicitness, we report only on specializing
 $f_{96}(j,x)$ to $j \in \scU_{3,1}(\Z[1/42])$.    From complete tables of elliptic curves \cite{CL07}, 
 this specialization set has size $413$.  Supporting Principle A, all $413$ algebras
 are non-isomorphic.  Supporting Principle B, these algebras all have Galois group 
 $A_{96}$.   Investigating Principle C is more subtle.  In lieu of completely factoring
 $f_{96}(j,x)$ over $\Q_p$ and taking field discriminants of the factors, we use 
 Newton polygons.    To illustrate this computationally much simpler method,
  we take $j = 1/3$ as a representative example, and work with 
 \[
 g(x) = 3 f(1/3,x) =  3^7 7^{21} x^{96} - 2^7 3^7 7^{21} x^{95} + \cdots + 2^{53} 3^{32} x  - 2^{48} 3^{31}.
 \]
 Factoring modulo $2$, $3$, and $7$ gives
 \begin{align*}
 g(x) & \stackrel{2}{\equiv} x^{96},
 &
 g(x) & \stackrel{3}{\equiv} x^{54} (x-2)^{33},
 &
 g(x) & \stackrel{7}{\equiv} h_{2}(x) h_{20}(x) h_{25}(x), 
 \end{align*}
 with $h_k(x)$ irreducible of degree $k$.  The $2$-adic Newton polygon of 
 $g(x)$ has all slopes $1/2$, showing that all $96$ roots $\alpha \in \overline{\Q}_2$ 
 have $\ord_2(\alpha)=1/2$.  Since the denominator is divisible by $2$, 
 one has that the $2$-adic wild degree as in \S\ref{wildram} is $m_{\rm 2 \mbox{-}wild} = 96$.   From a more
 complicated calculation with $3$-adic Newton polygons, we get that
 the $96$ roots $\alpha \in \overline{\Q}_3$ are distributed as follows:
 {\renewcommand{\arraycolsep}{1.7pt}
 \[
 \begin{array}{rrclrrcl}
 \mbox{9 roots with}  & \ord_3(\alpha) & = & 4/9, & \;\;\;  \mbox{22 roots with}  & \ord_3(\alpha-1) & = & 13/21, \\
 \mbox{27 roots with} & \ord_3(\alpha) & = &1/3, &  \mbox{3 roots with}  & \ord_3(3\alpha-(1-i)) & = & 2/3, \\
 \mbox{3 roots with} &  \ord_3(\alpha) & = & -1/3, &  \mbox{3 roots with}  & \ord_3( 3\alpha-(1+i)) & = & 2/3, \\
  \mbox{9 roots with} &  \ord_3(\alpha/3 - i) & = & 5/9, & \\
 \mbox{9 roots with} &  \ord_3(\alpha/3 + i) & = & 5/9, &  \mbox{12 roots with}  & \ord_3(\alpha-1) & = & 1/2. \\
 \end{array} 
 \]
 }
 
 \noindent Only the last twelve $\alpha$ could possibly not contribute to the $3$-adic wild degree,
 giving already $m_{\rm 3\mbox{-}wild} \geq 84$.  But in fact these $\alpha$ satisfy $\ord_3(\alpha^2/3-2)=5/12$ so one has $m_{\rm 3\mbox{-}wild} = 96$.   Finally, the $7$-adic Newton
polygon of $g(x)$ has slopes $0$ and $-3/7$ with multiplicities $47$ and $49$ respectively.
The slope of $0$ corresponds to the isolated roots modulo $7$ and the slope of $-3/7$ then
gives $m_{\rm 7\mbox{-}wild} = 49$.  

The Newton polygon process can be easily automated.  It says that all $413$ algebras are wildly ramified at both $2$ and $3$.  It says also that all algebras are wildly ramified at $7$ except
for those coming from the specialization points $-3^15^37^3/2^8$, $-7^3/2^13^2$,  
$7^3/2^9$, 
$7^3/3^5$, $7^3/2^13^3$, 
 $5^37^3/3^5$, $2^27^3/3$ and $7^4/2^63$, $-7^4/2^73^4$, $7^4$, $-7^5/2^13^8$.
 The first seven
all have tame ramification at $7$ corresponding to the partition $19^3 1^{39}$ while
the last four have tame ramification at $7$ corresponding to the partition
$57^1 3^{13}$.  This behavior comes from the fact that these specialization
points are all $7$-adically close
to $j=0$ and the degree $32$ polynomial above has tame ramification 
at $7$ given by the partition $19^1 1^{13}$.

  \section{A degree 202 family: degenerations and generic specialization} 
    \label{fam202}  \label{deg202}
          Continuing to increase degrees as we go through the last six
  sections, we now describe a family having degree $202$.    
  Our description emphasizes its degenerations, a relevant
  topic because how a family degenerates has substantial influence on
  how ramification behaves in the Hurwitz number fields within the family.  We conclude
  by observing that specialization is generic, both in one of the
  degenerations of the family and in the family itself.

  \subsection{Some plane curves}
  \label{planecurves}
  To streamline the subsequent subsections, we first present some polynomials
  defining affine curves 
    in the $x$-$y$ plane.  
  The next two subsections will place a natural function on each curve, 
  and we index the polynomials   by
  the degree of this function.     
  
   Eleven relatively simple polynomials are
  \begin{align*}
  A_{10}& =   x, & B_4 & = x-1, \\
  A_{13} & =  y, & B_8 & = y-1,  \\
  A_{14} & =  x-y, & B_{32} & = x^2 y-4 x^2-8 x y+20 x+10 y-20, \\
  A_{16} & =  3 x y-6 x-6 y+10, & C_{22} & = 3 x y^2-12 x y+8 x-15 y^2+40 y-24,\\
  A_{20} & =  x^2 y-3 x^2-6 x y  & C_{25} & = 3 x^2 y-6 x^2-12 x y+20 x+10 y-15, \\
              &   \qquad +12 x+6 y-10, & D_{10} & = 3 x^2-12 x+10.
  \end{align*}
  For all but one of these polynomials $P$, the curve $\AP$ given by its vanishing
  is obviously rational, as at least one of the variables appears to degree one
  in $P$.   In contrast, the curve $\AD_{10}$ consists of two genus
  zero components, neither one of which is defined over $\Q$.  
  
\begin{table}[htb]
  \[
   {\renewcommand{\arraycolsep}{2pt}
   \begin{array}{r|rrrrr}
   B_{52} &  1 & x & x^2 & x^3 & x^4  \\
   \hline
  1&                 1080 & -2160 & 1296 & -176 & -24 \\
   y &                 -1080 & 2052 & -1164 & 156 & 12 \\
  y^2 &                  135 & -276 & 180 & -36 &  \\
  y^3 &                 50 & -84 & 36 &  &  \\
  y^4 &                 15 & -12 &  &  &  \\
                  \end{array}
   \;\;\;\;\;\;\;\;\;           
    \begin{array}{r|rrrrr}
   D_{32} &  1 & x & x^2 & x^3 & x^4  \\
   \hline  
   1 &                 &  & 160 & -192 & 48 \\
   y &                 & -320 & 192 & 120 & -48 \\
   y^2 &                250 & 12 & -213 & 12 & 12 \\
   y^3 &               -300 & 288 & -18 & -12 &  \\
   y^4 &                 90 & -108 & 27 &  &  \\
   \end{array}
   }
  \]
  $\;$ \\
  \[
  \begin{array}{r|rrrrrrr}
 D_{48} &  1 & x & x^2 & x^3 & x^4 & x^5 & x^6 \\
 \hline
 1 &                   &  &  & 1600 & -2880 & 1632 & -288 \\
 y &                     &  & 2400 & -8160 & 9048 & -3960 & 576 \\
 y^2 &                    & 1200 & -6480 & 11448 & -8712 & 2880 & -324 \\
y^3 &                    -2500 & 6300 & -4620 & -108 & 1395 & -513 & 54 \\
y^4 &                   1500 & -3900 & 3780 & -1692 & 351 & -27 &  \\
                  \end{array}
  \]
  \caption{\label{bigpolys} Three polynomials $\sum c_{ij} x^i y^j$,
  presented by listing their coefficients $c_{ij}$.  }
  \end{table}
  
   Three more complicated polynomials are given in matrix form in Table~\ref{bigpolys}.
  The corresponding curves $\AB_{52}$, $\AD_{32}$, and $\AD_{48}$ have genus $1$, $2$, and $5$ respectively. 
  Each genus is much smaller than the upper bound allowed by the support in $\Z_{\geq 0}^2$ of the
  coefficients;  this bound, being the number of ``interior'' coefficients, is $6$, $6$, and $12$ respectively.  
  In each case, there are several singularities causing this genus reduction, one of which is
  at $(1,1)$.

  \subsection{Calculation of a rational presentation.}  This subsection is very similar to \S\ref{deg52comp}, illustrating that in favorable 
  cases computation of Hurwitz covers following the outline of \S\ref{computation} is quite mechanical.   As normalized Hurwitz 
  parameter we take  
 \[
 h = (S_6,(21111,3_021,3_1111,4_\infty11),(2,1,1,1)).
 \]
 Any function governed by $h$ is of the form
 \[
 g(s) = \frac{s^3 (s-x)^2 (s-y)}{a \left(s^2+s (d-e-1)+e\right)}.
 \]
  The ramification requirement on $g$ at $1$ is that $(g(1),g'(1),g''(1)) = (1,0,0)$.  
 These three equations let us express $a$, $d$, and $e$ in terms of 
$x$ and $y$.  Namely
  \begin{align*} 
  a & = -C_{25}, & d & = \frac{B^2_4 B_8}{C_{25}}, & e & = \frac{A_{20}}{C_{25}}. 
  \end{align*}
 Using a resolvent as usual, we find that the critical values
of $g(s)$ besides $0$, $1$, and $\infty$ are the roots of 
$W t^2+ (V-U-W) t + U t^2$ with 
\begin{eqnarray}
\label{202U} U & = & -2^2 3^3 A_{10}^5 A_{13}^4 A_{14}^3 A_{16} A_{20}^2,\\
\label{202V} V & = & -3^3 B_4 B_8^2 B_{32}^4 B_{52},  \\
\label{202W}W & = & 2^8 C_{25}^5 C_{22}.
\end{eqnarray}
Comparing with the standard quadratic $t^2 + (v-u-1) t + u$, 
one gets the rational presentation
\begin{align}
\label{ratpres202}
u & = \frac{U}{W}, & v & = \frac{V}{W}.
\end{align}
So, appealing to \eqref{202U}-\eqref{202W} and 
the explicit polynomials in \S\ref{planecurves},  Equations~\eqref{ratpres202} express 
$u$ and $v$ as explicit functions of $x$ and $y$.  

\subsection{A view of $\AX_h$}
\label{view}  Recall from \S\ref{specsixcurves} and Figure~\ref{spec2111}
that the complement of $\AU_{2,1,1,1}$ in the projective $u$-$v$-plane
consists of three lines $\AAnew$, $\AB$, $\AC$ and a conic $\AD$.  
   In the map from the affine $x$-$y$ plane to the projective $u$-$v$ plane, 
  we can consider the preimages of these discriminantal curves. 
  
     \begin{figure}[htb]
\includegraphics[width=4.6in]{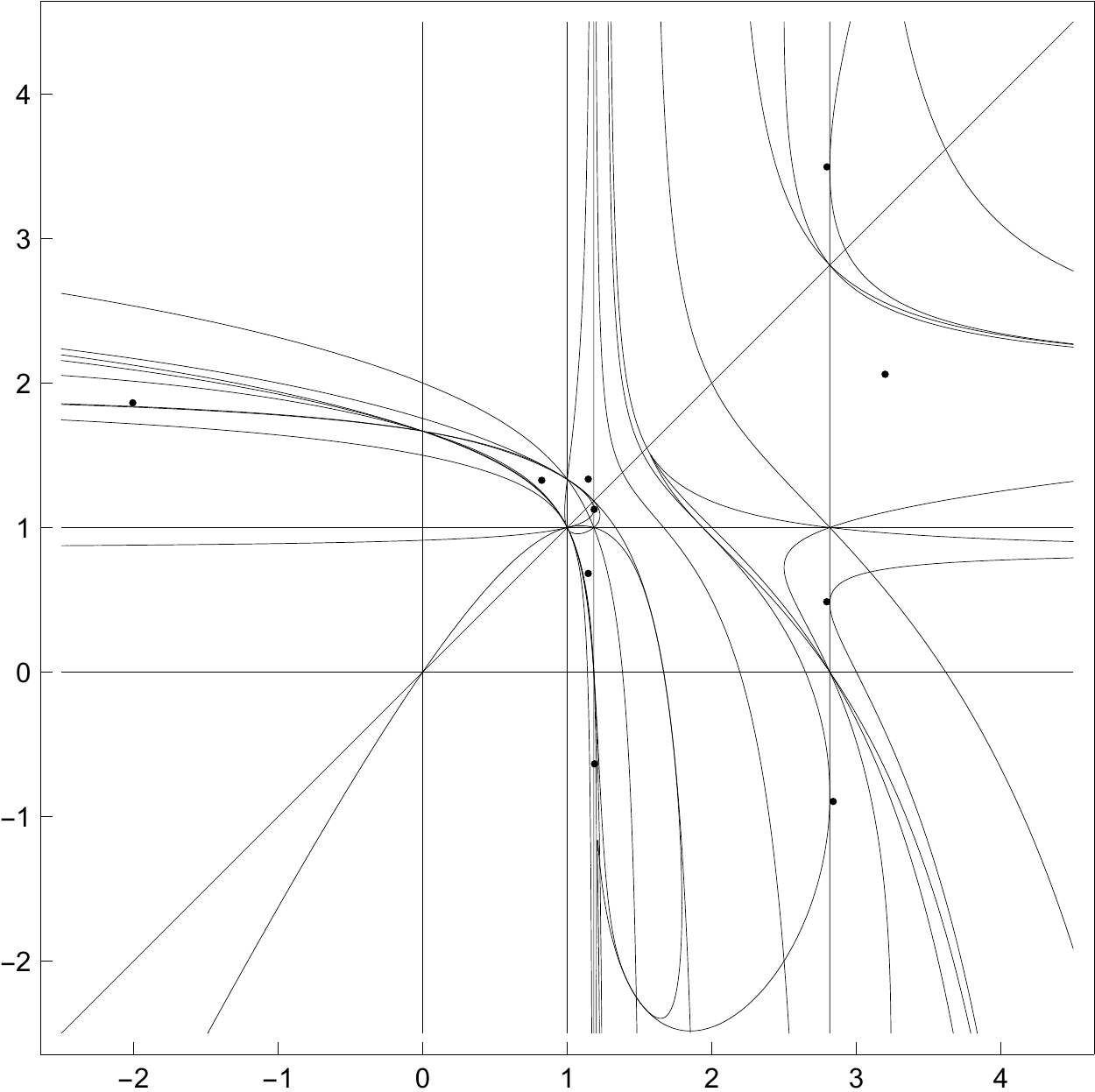}
\caption{\label{pict202}  $\scX_h(\R)$ is the complement of the drawn curves in the real $x$-$y$ plane. 
The drawn points are the ten real preimages of $(u,v) = (1,1)$.}
\end{figure}

    Figure~\ref{pict202} draws the real points of these four preimages. 
  Using as before a similar notation for an equation and its
  curve, inspection of our equations gives
  \begin{eqnarray*}
  \pi_h^{-1}(\AAnew) & = & \AAnew_{10} \cup \AAnew_{14} \cup \AAnew_{13}
   \cup \AAnew_{16} \cup \AAnew_{20}, \\
  \pi_h^{-1}(\AB) & = & \AB_4 \cup \AB_8 \cup \AB_{32} \cup \AB_{52}, \\
  \pi_h^{-1}(\AC) & = & \AC_{25} \cup \AC_{22},  \\
  \pi_h^{-1}(\AD) & = & \AD_{10} \cup \AD_{32} \cup \AD_{48}.
  \end{eqnarray*}
The figure is intended to indicate the rich geometry present in any Hurwitz
surface.  Other interesting curves present whenever $\nu = (2,1,1,1)$ are the 
preimages of the lines ${\sf ad}$, ${\sf bd}$, ${\sf cd}$, ${\sf bc}$, ${\sf ac}$, and ${\sf ab}$ 
introduced
in \S\ref{deg52}.  For the current $h$, all of them have a complicated 
real locus.  Their genera are respectively $25$, $18$, $23$, $35$, $31$, and $23$.
The curves ${\sf ad}$, ${\sf bd}$, and ${\sf cd}$ intersect at the preimage of the point ${\sf d}$, and 
Figure~\ref{pict202} also draws the ten real points of this preimage. 

\subsection{Degree $202$ polynomials and their degenerate factorizations}  Removing $y$ and $x$ respectively 
from \eqref{ratpres202} by resultants gives degree $202$ polynomials
$f(u,v,x)$ and $\phi(u,v,y)$.   Completely expanded, they 
have 10484 and 15555 terms respectively.  

The structures studied in the previous two subsections appear when one factors
specializations corresponding to the four discriminantal components:
\begin{eqnarray*}
f(0,v,x) & = & x^{50} (x^2-4x+6) a_{13}(v,x)^4 a_{14}(v,x)^3 a_{16}(v,x) a_{20}(v,x)^2, \\
f(u,0,x) & = & -(x-1)^6 b_8(u,x)^2 b_{32}(u,x)^4 b_{52}(u,x), \\
\lim_{u \rightarrow \infty} \frac{f(u,w u,x)}{u^{10}} & = & -2^{10} (x-1)^7 c_{10}(w,x)^3 c_{22}(w,x) c_{25}(w,x)^5, \\
f(r^2,(1-r)^2,x) & = & -(3x^2-12 x+10)^5 d_{32}(r,x)^3 d_{48}(r,x)^2,
\end{eqnarray*}
and
\begin{eqnarray*}
\phi(0,v,y) & = & y^{52} (3 y^2 - 8 y + 8) \alpha_{10}(v,y)^5  \alpha_{14}(v,y)^3 \alpha_{16}(v,y) \alpha_{20}(v,y)^2, \\
\phi(u,0,y) & = & (y-1)^{16} (y-4)^2  \beta_4(u,y)\beta_{32}(u,y)^4 \beta_{52}(u,y), \\
\lim_{u \rightarrow \infty} \frac{\phi(u,w u,y)}{u^{13}} & = & -2^{36} (y-2)^3  (y-1)^7 \gamma_{6}(w,y)^2 \gamma_{22}(w,y) \gamma_{25}(w,y)^5, \\
\phi(r^2,(1-r)^2,y) & = & \delta_{10}(r,y) \delta_{32}(r,y)^3 \delta_{48}(r,y)^2.
\end{eqnarray*}
Our notation coordinates the different viewpoints: for example, the equations $D_{48}=0$,
$d_{48}(r,x)=0$, and $\delta_{48}(r,y)=0$ all describe the genus five curve ${\AD}_{48}$.

  As a sample degeneration, chosen because it makes an interesting
   comparison the degree $25$ polynomials from our introductory example, 
 \begin{eqnarray*}
 \lefteqn{c_{25}(w,x)=} \\
&& 
-(2 x-5) (3 x-5)^2 \left(6 x^4-40 x^3+105 x^2-120 x+50\right)^4  \\
&& \qquad  \left(12 x^6-60 x^5-40 x^4+760 x^3-1800 x^2+1750 x-625\right) \\
&& + 4 w x^5 (x^2-5 x+5)^3 (3 x^2-10 x+10)^2 (6 x^2-20
    x+15)^4 (6 x^2-15 x+10).
\end{eqnarray*}
Here the $x$-line is identified with $\AX_h$ for 
$h = (S_6,(321,3111,51,21111),(1,1,1,1))$, and
the the $w$-line with 
$\AU_{1,1,1,1}$.   All the other degenerations have a similar four-point
description.    The discriminant of $c_{25}(w,x)$ is $2^{212} 3^{66} 5^{285} w^{13} (w-1)^{19}$.  
All the other degenerations are likewise three-point covers,  all full except for 
$\AAnew_{14}$, $\AAnew_{16}$,
$\AAnew_{20}$, $\AC_6$, and $\AD_{10}$.     
In every case, the target variable, be it $v$, $u$, $w$, or $r$, is chosen
such that the singular values are $0$, $1$, and $\infty$.  

  To be noted is that we are not expending any extra effort here to introduce
 a conceptually defined completion of $\AX_h$.  Indeed the curves
 that consist of horizontal lines, namely $\AAnew_{13}$ and $\AB_{8}$, are 
 seen clearly by $f$ but only as vestigial factors by $\phi$.  In 
 reverse, the curves that consist of vertical lines, namely $\AAnew_{10}$, $\AB_4$, and 
 $\AD_{10}$, 
 are seen completely by $\phi$ but only partially by $f$.  Finally,
 to see preimages corresponding to the factors $c_{10}(w,x)$ and
 $\gamma_6(w,x)$, one would have to go beyond the $x$-$y$ plane
 as a partial completion of $\AX_h$.

    A braid group computation gives the partition of 202 which 
 captures how local sheets of $\AX_h$ are interchanged as
 one goes around one of the four discriminantal components
 in the completion of $\AU_{2,1,1,1}$.   These partitions are
 \begin{align*}
 \beta_A & = 5^{10} 4^{13} 3^{14} 2^{20} 1^{16 + \mathbf{2}}, &
 \beta_B & =  4^{32} 2^{8+\mathbf{1}} 1^{52 + 4}, &
 \beta_C & =  5^{25} 3^{10+\mathbf{2}} 2^{6+\mathbf{2}} 1^{22+ \mathbf{3}}, 
  \end{align*}
 and $\beta_D  =  3^{32} 2^{48} 1^{10}$.  The boldface
 exponents correspond to components not seen by our simple 
 calculations.  Thus we are missing only $2$, $2$, $13$, and $0$ 
 of the $202$ sheets near the preimages of $\AAnew$, $\AB$, $\AC$, and $\AD$ respectively.

 \subsection{Specialization}  The degenerations can be specialized, and the computations support Principles A, B, and C.  For example, consider
 $c_{25}(w,x)$ specialized to $w$ in the known set $\scU_{1,1,1,1}(\Z[1/30])$.  The $99$ algebras 
 are all distinct, they are all full, and they are all wildly ramified at each of $2$, $3$, and 
 $5$.   The wild ramification is even less than in our introductory example
 described in and around Figure~\ref{discfchart} for two reasons.  
 First, the algebras decompose more here;
 in particular, there are no $2$-adic factors of degree $16$ and no $5$-adic factors
 of degree $25$.  Second, the individual factors are less ramified; for example,
 $2$-adic octics have a largest discriminant of $2^{26}$ rather than $2^{31}$
 and $3$-adic nonics have a largest discriminant of $3^{16}$ rather than $3^{26}$.  
 The mean exponents are now 
 $(\langle a \rangle, \langle b \rangle, \langle c \rangle) \approx (47,32,35)$ 
 rather than $(56,43,42)$; the maximum exponents are now 
 $(a_{\rm max},b_{\rm max},c_{\rm max}) = (61,38,44)$ rather than
 $(79,57,52)$.
 
 Specialization of the full family at the $2947$ points of $\scU_{2,1,1,1}(\Z[1/30])$ 
 can also be satisfactorily studied, despite 
 the large degree.   The $2947$ algebras are all distinct and they all have 
 Galois group $A_{202}$ or
  $S_{202}$.  From Newton polygons, we know they are all wildly ramified 
  at $2$, $3$, and $5$.  {  Thus, in this family, Principles A, B, and C hold without
  exception. }

\section{A degree 1200 field: 
computations in large degree } 
\label{deg1200} \label{largedegree}
      Conjecture~\ref{mc} says that for certain finite sets of primes $\cP$, there exist full number fields
of arbitrarily large degree with ramification set in $\cP$.   A natural computational challenge
for a given $\cP$ is then to produce an explicit full Hurwitz number field $K_{h,u}$ 
with degree $m$ as large as possible.   In this short final section,
we take $\cP=\{2,3,5\}$ and produce such a field for degree $m=1200$.  

Taking 
$h = (S_6, (21111,3_02_11, 4_\infty11),(4,1,1))$ and normalizing as indicated,
the functions to consider are 
\[
g(s) = \frac{a s^3 (s-1)^2 (s-x)}{s^2 + b s+c}.
\]
As specialization point, we take $u = (( t^4 -4 t -6),\{0\},\{\infty\})$.
 This specialization point
indeed keeps ramification within $\{2,3,5\}$ as the
discriminant of $ t^4 -4 t -6$ is
$-2^8 3^5$.  

The condition that the critical values besides $0$ and $\infty$ are the roots of $t^4 -4 t -6$ gives
 four equations in the four unknowns $x$, $a$, $b$, $c$.  
Of the unknowns, we are focusing on $x$ because its special values $0$,
$1$, and $\infty$ are all meaningful, corresponding to degenerations.  
Eliminating $a$ and then $c$ is easy.  Eliminating $b$ then
has a ten-minute run-time on {\em Magma} to get a degree 3700 polynomial.  
Factorizing this polynomial to find the relevant factor has a 
one-minute run-time.     The resulting monic polynomial
$f_{1200}(x) \in \Z[1/30][x]$ defining $K_{h,u}$ satisfies $f_{1200}(0) = 2^{880}/5^{500}$
and $f_{1200}(1) = 3^{684}/2^{256} 5^{500}$.   
After removing all factors of $2$, $3$, and $5$, 
the coefficients are integers averaging about $440$ digits. 

\begin{figure}[htb]
\includegraphics[width=4.6in]{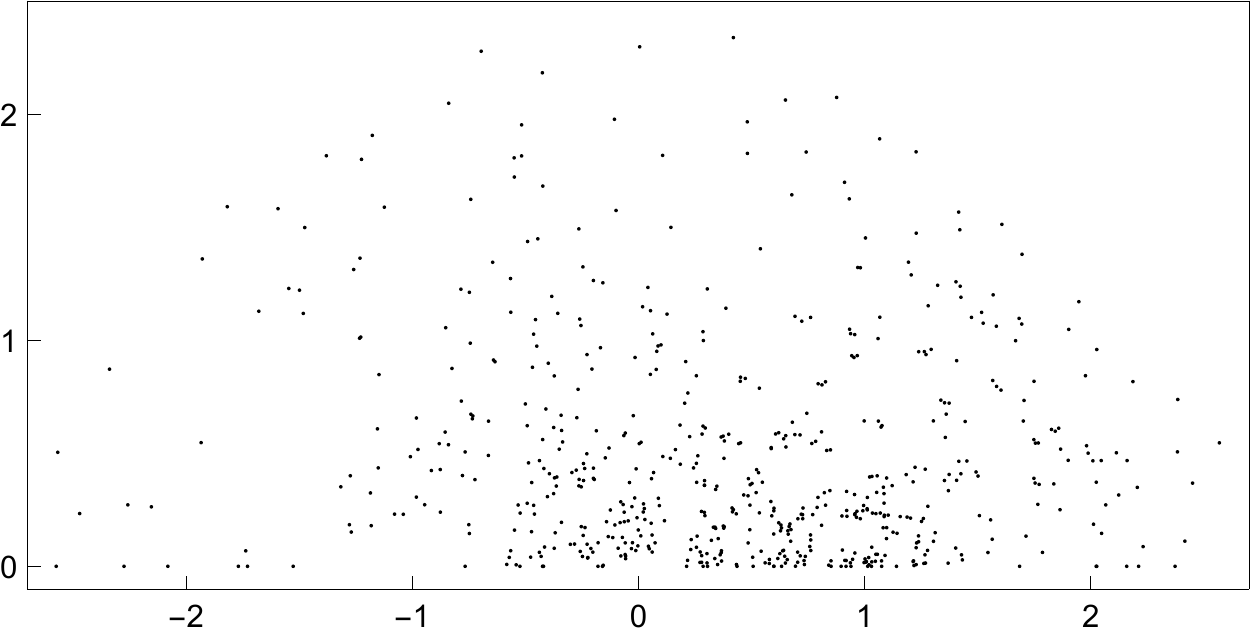}
\caption{\label{pict1200}  Roots in the closed upper half plane of a polynomial defining a degree $1200$ Hurwitz number field}
\end{figure}

The polynomial is to some extent analyzable despite its large degree and large coefficients.  
From the factorization partitions
$(989, 208, 3)$ at $19$ and 
$(1181,9,6,4)$ at $47$, it has Galois group $S_{1200}$,
in conformity with Principle B.  From Newton polygons, it is wildly
ramified at $2$, $3$, and $5$, as predicted by Principle C.  
Figure \ref{pict1200} presents the roots of $f_{1200}(x)$ in the
closed upper half plane, all of which lie in the drawn
window $[-2.6,2.6] \times [0,2.4]$.  There 
are $34$ real roots and a general tendency
of roots to cluster near the interval $[0,1]$ connecting
the special points $0$ and $1$.  

Large degree Hurwitz number fields 
provide specific challenges 
to improve computational algorithms for general number fields.
For example, from Newton polygons we have substantial information on how the
field $K_{h,u}$ of this section factors over $\Q_2$, $\Q_3$, and $\Q_5$.  
However we cannot go far enough to determine the exponents in
its discriminant $-2^a 3^b 5^c$.

 \bibliographystyle{amsplain}
\bibliography{mast3}

\providecommand{\bysame}{\leavevmode\hbox to3em{\hrulefill}\thinspace}
\providecommand{\MR}{\relax\ifhmode\unskip\space\fi MR }
\providecommand{\MRhref}[2]{%
  \href{http://www.ams.org/mathscinet-getitem?mr=#1}{#2}
}
\providecommand{\href}[2]{#2}
\begin{thebibliography}{10}

\bibitem{BR11}
Jos{\'e} Bertin and Matthieu Romagny, \emph{Champs de {H}urwitz}, M\'em. Soc.
  Math. Fr. (N.S.) (2011), no.~125-126, 219. \MR{2920693}

\bibitem{Bha05}
Manjul Bhargava, \emph{The density of discriminants of quartic rings and
  fields}, Ann. of Math. (2) \textbf{162} (2005), no.~2, 1031--1063.
  \MR{2183288}

\bibitem{Bha07}
\bysame, \emph{Mass formulae for extensions of local fields, and conjectures on
  the density of number field discriminants}, Int. Math. Res. Not. IMRN (2007),
  no.~17, Art. ID rnm052, 20. \MR{2354798}

\bibitem{Bha10}
\bysame, \emph{The density of discriminants of quintic rings and fields.}, Ann.
  of Math. (2) \textbf{172} (2010), no.~3, 1559--1591.

\bibitem{magma}
Wieb Bosma, John Cannon, and Catherine Playoust, \emph{The {M}agma algebra
  system. {I}. {T}he user language}, J. Symbolic Comput. \textbf{24} (1997),
  no.~3-4, 235--265, Computational algebra and number theory (London, 1993).
  \MR{MR1484478}

\bibitem{CL07}
J.~E. Cremona and M.~P. Lingham, \emph{Finding all elliptic curves with good
  reduction outside a given set of primes}, Experiment. Math. \textbf{16}
  (2007), no.~3, 303--312. \MR{2367320}

\bibitem{DH71}
H.~Davenport and H.~Heilbronn, \emph{On the density of discriminants of cubic
  fields. {II}}, Proc. Roy. Soc. London Ser. A \textbf{322} (1971), no.~1551,
  405--420. \MR{0491593}

\bibitem{FK14}
Claus Fieker and J{\"u}rgen Kl{\"u}ners, \emph{Computation of {G}alois groups
  of rational polynomials}, LMS J. Comput. Math. \textbf{17} (2014), no.~1,
  141--158. \MR{3230862}

\bibitem{JR14}
John~W. Jones and David~P. Roberts, \emph{A database of number fields}, LMS J.
  Comput. Math. \textbf{17} (2014), no.~1, 595--618, Database at
  \url{http://hobbes.la.asu.edu/NFDB/}. \MR{3356048}

\bibitem{Kat96}
Nicholas~M. Katz, \emph{Rigid {L}ocal {S}ystems}, Annals of Mathematics
  Studies, vol. 139, Princeton University Press, Princeton, NJ, 1996.
  \MR{1366651}

\bibitem{Kon}
Joachim Koenig, \emph{Computation of {H}urwitz spaces and new explicit
  polynomials for almost simple {G}alois groups}, Preprint, to appear in Math.\
  Comp.\, 2015.

\bibitem{Mal00}
Gunter Malle, \emph{Multi-parameter polynomials with given {G}alois group}, J.
  Symbolic Comput. \textbf{30} (2000), no.~6, 717--731. \MR{1800034}

\bibitem{MM99}
Gunter Malle and B.~Heinrich Matzat, \emph{Inverse {G}alois theory}, Springer
  Monographs in Mathematics, Springer-Verlag, Berlin, 1999. \MR{1711577}

\bibitem{MR05}
Gunter Malle and David~P. Roberts, \emph{Number fields with discriminant {$\pm
  2^a3^b$} and {G}alois group {$A_n$} or {$S_n$}}, LMS J. Comput. Math.
  \textbf{8} (2005), 80--101 (electronic). \MR{2135031}

\bibitem{RobHBM}
David~P. Roberts, \emph{Hurwitz-{B}elyi maps}, Arxiv, August 30, 2016.
  Submitted.

\bibitem{Rob07}
\bysame, \emph{Wild partitions and number theory}, J. Integer Seq. \textbf{10}
  (2007), no.~6, Article 07.6.6, 34. \MR{2335791}

\bibitem{RobG2}
\bysame, \emph{Division polynomials with {G}alois group {$SU_3(3).2\cong
  G_2(2)$}}, Advances in the theory of numbers, Fields Inst. Commun., vol.~77,
  Fields Inst. Res. Math. Sci., Toronto, ON, 2015, pp.~169--206. \MR{3409329}

\bibitem{Rob15}
\bysame, \emph{Polynomials with prescribed bad primes}, Int. J. Number Theory
  \textbf{11} (2015), no.~4, 1115--1148. \MR{3340686}

\bibitem{RV15}
David~P. Roberts and Akshay Venkatesh, \emph{Hurwitz monodromy and full number
  fields}, Algebra Number Theory \textbf{9} (2015), no.~3, 511--545.
  \MR{3340543}

\bibitem{Se08}
Jean-Pierre Serre, \emph{Topics in {G}alois theory}, second ed., Research Notes
  in Mathematics, vol.~1, A K Peters, Ltd., Wellesley, MA, 2008, With notes by
  Henri Darmon. \MR{2363329}

\bibitem{Pari}
{The PARI~Group}, Bordeaux, \emph{{PARI/GP version {\tt 2.7.5}}}, 2015,
  available from \url{http://pari.math.u-bordeaux.fr/}.

\bibitem{Vol01}
Helmut V{\"o}lklein, \emph{A transformation principle for covers of {${\Bbb
  P}^1$}}, J. Reine Angew. Math. \textbf{534} (2001), 156--168. \MR{1831635}

\bibitem{math}
{Wolfram Research, Inc.}, \emph{Mathematica 10.0.2.0}.

\end{thebibliography}

\end{document}